\documentclass[a4paper,12pt]{article}

\usepackage{hyperref}
\usepackage[vmargin=2cm,hmargin=2cm,headheight=14.5pt,top=2cm,headsep=.5cm]{geometry}

\usepackage{bm}
\usepackage[utf8]{inputenc}

\usepackage[autobold]{mathfixs}
\usepackage{empheq}
\usepackage{stackrel}
\usepackage{cases}
\usepackage{mathtools}
\usepackage{amsthm,amsmath,amscd}
\usepackage[charter]{mathdesign}
\usepackage{perpage}
\usepackage{bm}
\usepackage{tikz-cd}
\tikzcdset{every label/.append style = {font = \small}}
\usepackage{caption}
\usepackage{graphicx}
\DeclareMathSizes{12}{12}{8}{6}

\usepackage{cite}
\usepackage{enumitem}
\usepackage{framed}
\usepackage[symbol]{footmisc}
\usepackage{marvosym}
\usepackage{perpage}
\MakePerPage[1]{footnote}

\makeatletter
\renewcommand*{\@makefnmark}{\hbox{\@textsuperscript{%
			\normalfont\@thefnmark}}}
\makeatother
\makeatletter
\def\@fnsymbol#1{\ensuremath{\ifcase#1\or \text{\Mercury}
		\or \text{\Venus} \or \text{\Earth} \or
		\text{\Jupiter} \or \text{\Saturn} \or \text{\Neptune} \or \text{\Uranus} \or
		\text{\Pluto}
		\or \text{\Moon} \or \text{\Sun}
		\else\@ctrerr\fi}}%
\makeatother

\newtheoremstyle{ptheorem}{1em}{0em}{\itshape}{}{\bfseries}{.}{.5em}{\thmname{#1}\thmnumber{
		#2}\thmnote{ (\hspace{-.01pt}{#3})}}

\theoremstyle{ptheorem}

\newtheorem{thm}{Theorem}[section]
\newtheorem{pro}[thm]{Proposition}
\newtheorem{lem}[thm]{Lemma}

\newtheoremstyle{hdef}{1em}{0em}{}{}{\bfseries}{.}{.5em}{\thmname{#1}\thmnumber{
		#2}\thmnote{ (\hspace{-.01pt}{#3})}}
\theoremstyle{hdef}

\newtheorem{dfn}[thm]{Definition}
\newtheorem{rem}[thm]{Remark}

\newtheorem{exa}[thm]{Example}

\numberwithin{equation}{section}
\numberwithin{figure}{section}

\let\originalleft\left
\let\originalright\right
\renewcommand{\left}{\mathopen{}\mathclose\bgroup\originalleft}
\renewcommand{\right}{\aftergroup\egroup\originalright}

\DeclareMathOperator{\supp}{supp}

\DeclareMathOperator{\dif}{d}

\newcommand{\cB}{{\mathcal B}}
\newcommand{\cC}{{\mathcal C}}

\newcommand{\cP}{{\mathcal P}}

\newcommand{\cT}{{\mathcal T}}
\newcommand{\cU}{{\mathcal U}}
\newcommand{\cV}{{\mathcal V}}

\newcommand{\bN}{{\mathbb N}}

\newcommand{\bP}{{\mathbb P}}

\newcommand{\bR}{{\mathbb R}}
\newcommand{\bS}{{\mathbb S}}

\renewcommand{\a}{\alpha}
\renewcommand{\b}{\beta}
\renewcommand{\c}{\gamma}
\newcommand{\C}{\Gamma}

\newcommand{\e}{\varepsilon}

\renewcommand{\k}{\kappa}

\renewcommand{\phi}{\varphi}
\renewcommand{\le}{\leqslant}
\renewcommand{\ge}{\geqslant}

\newcommand{\ol}{\overline}
\newcommand{\fa}{\forall}

\newcommand{\n}{{n\in\bN}}

\newcommand{\Ra}{\Rightarrow}

\newcommand{\nkp}{\enskip}

\newcommand{\sfa}{\nkp\fa}
\renewcommand{\d}{\delta}

\renewcommand{\(}{\left(}
\renewcommand{\)}{\right)}

\newcommand{\lil}{\lim\limits}
\newcommand{\til}{\widetilde}
\newcommand{\pd}{\partial}

\newcommand{\bs}{\backslash}
\newcommand{\Lsp}[1]{\operatorname{L^{#1}}}
\newcommand{\olb}[1]{%
	\vbox{\offinterlineskip\ialign{\hfil##\hfil\cr
			$\rotatebox[origin=c]{90}{$]$}$\cr\noalign{\kern-.45ex}{$#1$}\cr}}}

\newcommand{\noop}[1]{}

\renewcommand{\ss}{\subset}

\allowdisplaybreaks

\title{Asymptotic properties of PDEs in compact spaces}

\date{}
\author{Lucía López-Somoza and F. Adri\'an F. Tojo}

\begin{document}
	\maketitle

	\begin{center}  {\itshape\small Departamento de Estat\'{\i}stica, An\'alise Matem\'atica e Optimizaci\'on \\ Instituto de Matem\'aticas \\ Universidade de Santiago de Compostela \\ 15782, Facultade de Matem\'aticas, Campus Vida, Santiago, Spain. \\ e-mail: lucia.lopez.somoza@usc.es, fernandoadrian.fernandez@usc.es}
	\end{center}

	\medbreak

\begin{abstract}
		In this article we combine  the study of solutions of PDEs with the study of asymptotic properties of the solutions via compactification of the domain. We define new spaces of functions on which study the equations, prove a version of Ascoli-Arzelà Theorem, develop the fixed point index results necessary to prove existence and multiplicity of solutions in these spaces and also illustrate the applicability of the theory with an example.
\end{abstract}

	\medbreak


	\medbreak


\section{Introduction}

The use of topological methods in the study of PDE is a classical field of research --see \cite{Gilbarg,GeorgZen,GrossinhoTersian,kawanago}. Unfortunately, the most modern and sophisticated methods that have been recently developed for ODEs --see, for instance, \cite{Bai2020,CabadaDimitrov2021,Infante2021a,Sousa2021}-- have been difficult to apply to PDEs. The reasons for this are, to cite some, the greater effort needed to check the, if rather weak, cumbersome hypotheses, the lower availability of explicit expressions of Green's functions for PDEs, the higher complexity of the domain of definition and the higher regularity that is necessary in order to obtain existence and uniqueness results.

Even then, there has been a recent effort to overcome this difficulties, mainly by imposing some kind of symmetry on the operator that defines the equation and, in particular, searching for radial solutions \cite{Goodrich2021,Goodrich2018,Goodrich2018a,Infante2015,Cianciaruso2019,Cianciaruso2018,Goodrich2017,Goodrich2016}. More general approaches also appear for elliptic PDEs and systems of PDEs \cite{Infante2021,Infante2018,Infante2019}.

On the other hand, the study of ODEs on unbounded domains has progressed steadily \cite{MinCar,MinCar2,MinCar-4ord,FiMinCar,DjeGue}. The key to deal with unbounded domains is to use some kind of relatively compactness criterion such as \cite[Theorem 1]{przerad} --see for instance \cite{Corduneanu,przerad}. These kind of criteria are reworkings of the classical Ascoli-Arzelà Theorem and have been used in a different way in \cite{Cabada2018,Cabada2020}. In these works the authors are able to apply Ascoli-Arzelà Theorem by compactifying the domain of the functions involved in the ODE, thus allowing for a study of the asymptotic properties of the solutions.

In this article we combine both the study of solutions of PDEs with the study of asymptotic properties of the solutions via compactification of the domain. Furthermore, we take the opportunity to fix some of the shortcomings in \cite{Cabada2018,Cabada2020} and provide an example of application. It is worth noticing that our results are not constrained to partial differential equations of a particular type (parabolic, elliptic, hyperbolic) since the results obtained are presented for the  integral form of the equations; but also that, in general, the conditions to be checked for a particular problem can become quite unwieldy, which can be a limiting factor when it comes to apply the results to more convoluted problems.

The structure of this article is as follows. On Section 2 we deal with the basic topological notions necessary for understanding compactifications and provide some examples thereof. In Section 3 we provide the definition of the family of Banach spaces we will be dealing with. We also prove basic results regarding its structure as Banach space as well as a version of Ascoli-Arzelà Theorem (Theorem~\ref{taae}) necessary for the results to come. It is in Section 4 that we apply the usual topological methods regarding the fixed point index in order to obtain existence results of an integral problem in several variables which, in general, can be seen as a transformation on a PDE problem. Finally, in Section 5 we provide an example of a hyperbolic equation of which a solution with a predetermined asymptotic behavior can be found.

\section{Preliminaries: compactifications and extensions}


In order to understand asymptotic behavior on a metric space $(X,d)$ we first need to formalize the notion of \emph{point of infinity}. In an intuitive way, we can picture a point of infinity as a point far away from any  point in $X$. For instance, the usual order relation on the real numbers makes us think of a number bigger than any other. Those points of infinity must live some place, and that place is what we call a \emph{compactification}. A compactification $\til X$ has a topological structure that allows us to formalize the notion of asymptotic behavior of those functions defined on $X$ in a precise way. This is because, once we have a topology in $\til X$, we can take limits. Furthermore, the topology of the compactification, when metrizable, allows us to study the relations and relative positions of the different points of infinity.

We proceed now to formally define the concept of compactification through an adequate map.

\begin{dfn}Let $X$, $Y$ be topological spaces, $Y$ compact. We say that a continuous function $\kappa:X\to Y$ is a \emph{compactification} of $X$ if $\kappa(X)$ is dense in $Y$ (that is, $\ol{\kappa(X)}=Y$) and $\kappa: X\to\kappa(X)$ is a homeomorphism. We will usually identify the compactification with $Y$. Also, we will denote the inverse of $\k:X\to\k(X)$ as $\k^{-1}|_{\k(X)}$.
\end{dfn}

In this work we will restrict ourselves to metric compactifications, that is, to the case where $Y$ is a compact metric space. The interested reader may find more information regarding compactifications and their properties in \cite{Willard,lynn,Dugundji,chandler,Kura}.
\begin{rem} In practice, any compactification $\kappa:X\to Y$ can be considered to be the identitary inclusion. To see this, assume $\tau$ is the topology of $Y$ and consider the bijective map
	\begin{center}\begin{tikzcd}[row sep=tiny]
\cP(X\sqcup (Y\bs \kappa(X)))\arrow{r}{\Theta} &  \cP(Y)\\
U\sqcup V\arrow[hook]{r} & \kappa(U)\cup V
	\end{tikzcd}
	\end{center}
where $\sqcup$ denotes the topological sum and $\cP(\cdot)$ is the power set. $\Theta^{-1}(\tau)$ is a topology in 	$X\sqcup (Y\bs \kappa(X))$ that makes $Y$ and $X\sqcup (Y\bs \kappa(X))$ homeomorphic with the homeomorphism $\xi:X\sqcup (Y\bs \kappa(X))\to Y$ defined as
\[\xi(x)=\begin{dcases} \kappa(x), & x\in X,\\ x, & x\in Y\bs \kappa(X).\end{dcases}\] Furthermore, $\til\kappa:X\to X\sqcup (Y\bs \kappa(X))$ such that $\til\kappa(x)=x$ is a compactification and $\kappa=\xi\circ\til\kappa$.
\end{rem}

In order to illustrate the notion of compactification we present now various examples where the compactifications chosen map the space $\bR^n$ to well known differentiable manifolds. It is  important to point out that any compact differentiable manifold is a compactification of $\bR^n$ as a consequence of\cite[Corollary~2.8, p. 271]{doCarmo1992}, so this is a very general situation.

\begin{exa}[Directional compactification]Let $n\in\bN$, $B:=\{x\in\bR^n\ :\ \|x\|\le1\}$. $\kappa:\bR^n\to B$ defined as $\kappa(x)=x/(1+\|x\|)$. $\kappa:\bR^n\to \mathring B$ is a $\cC^\infty$-diffeomorphism and $\ol{\kappa(\bR^n)}=B$, so $\kappa$ is a compactification.

    Observe that the elements of $\partial B=\bS^{n-1}$ denote different \emph{`directional points of infinity'} in the sense that, if $v\in\bS^{n-1}$ and $f(t)=tv$, $t\in\bR$, then $\lim_{t\to\infty}\kappa(f(t))=v$.
\end{exa}
\begin{exa}[Projective spaces] Take $\kappa:\bR^n\to B$ as before. Now, we establish an equivalence class in $B$ in the following way: $x\sim y$ iff $x=y$ or $x,y\in\partial B$ and $x=-y$. With this equivalence class $B|_\sim$ is homeomorphic to the $n$-th real projective space $\bP^n$.  Let $\pi:B\to B|_\sim$ be the projection onto the quotient space. We can consider the compactification $\pi\circ\kappa:\bR^n\to B|_\sim$.
\end{exa}
\begin{exa}[Alexandroff's one-point compactification] Again, take $\kappa:\bR^n\to B$ as before and consider the equivalence class in $B$ defined by $x\sim y$ iff $x=y$ or $x,y\in\partial B$. With this equivalence class $\partial B$ reduces to a point and $B|_\sim$ is homeomorphic to the $n$-th sphere $\bS^n$. If $\pi:B\to B|_\sim$ is the projection onto the quotient space, we can consider the compactification $\pi\circ\kappa:\bR^n\to B|_\sim$.
\end{exa}

    The key to develop our theory and relate it to the concepts ahead is an adequate notion of limit which we present now.

\begin{dfn}\label{dfnlim}Let $X$, $Y$, $Z$ be topological spaces, $Y$ compact, $Z$ Hausdorff,  $\kappa:X\to Y$ a compactification of $X$, $y\in Y\bs \kappa(X)$ and $f:X\to Z$. We say the \emph{limit of $f$ when $x$ tends to $y$ is $z\in Z$}, and we write $\lim_{x\to y}^\kappa f(x)=z$, if for every neighborhood $V$ of $z$ in $Z$ there exists a neighborhood $U$ of $y$ in $Y$ such that $f(\kappa^{-1}(U\bs\{y\}))\subset V$.
\end{dfn}
\begin{rem} Since $y\in Y\bs \kappa(X)$, we have that $\kappa^{-1}(U\bs\{y\})=\kappa^{-1}(U)$, so we could have written $f(\kappa^{-1}(U))\subset V$ in Definition~\ref{dfnlim}.
\end{rem}
Observe that Definition~\ref{dfnlim} is dependent on the compactification $\kappa$ as the following example illustrates.
\begin{exa}
	Let $X=[-\infty,\infty]$ with the usual compact interval topology and $Y=X|_{\pm\infty}$, that is, the quotient of $X$ by the relation that identifies $-\infty$ and $\infty$. We can consider the following two compactifications of $\bR$, $\kappa_1:\bR\to X$ and $\kappa_2:\bR\to Y$ given by $\kappa_1(t)=\kappa_2(t)=t$. Then, if we consider $f(x)=\arctan x$, we have that $\lim_{x\to \infty}^{\kappa_1} f(x)=\lim_{x\to \infty} f(x)=\frac{\pi}{2}$ but $\lim_{x\to \infty}^{\kappa_2} f(x)=\lim_{x\to \pm\infty} f(x)$ does not exist.
\end{exa}

Definition~\ref{dfnlim} has a simpler form in the case of metric spaces. In the following we will denote by $d$ the distance in any metric space.
\begin{pro}\label{prolim}Let $X$, $Y$ and $Z$ be metric spaces, $Y$ compact,  $\kappa:X\to Y$ a compactification of $X$, $y\in Y\bs \kappa(X)$ and $f:X\to Z$. Then $\lim_{x\to y}^\kappa f(x)=z$ if and only if for every $\e\in\bR^+$ there exists $\d\in\bR^+$ such that $d(f(x),z)<\e$ if $d(\kappa(x),y)<\d$.
\end{pro}

 We will be using the following result.
\begin{pro}[{\cite[Proposition 1.29]{chandler}}]\label{proun} Let $X$, $Y$ be topological spaces, $A\subset X$, $f:A\to Y$ continuous. If $\til f:X\to Y$ is a continuous extension of $f$, $A$ is dense in $X$ and $Y$ is Hausdorff, then $\til f$ is unique.
\end{pro}


\begin{thm}[Existence of continuous extensions]\label{lemex}Let $X$, $Y$ and $Z$ be metric spaces, $Y$ compact, $\kappa:X\to Y$ a compactification of $X$ and $f:X\to Z$ continuous. Then the following statements are equivalent:
\begin{enumerate}
\item There exists $\lim_{x\to y}^\kappa f(x)$ for every $y\in Y\bs \kappa(X)$.
\item There exists a continuous map $\til f:Y\to Z$ such that $f=\til f\circ\kappa$.
\end{enumerate}
Furthermore, the extension $\til f$ is unique.
\end{thm}
\begin{proof}

    (I)$\Ra$(II) Since $Z$ is Hausdorff, the limit $\lim_{x\to y}^\kappa f(x)$ is unique, so we can define
    \[\til f(y):=\begin{dcases} f(\kappa^{-1}(y)), & y\in \kappa(X),\\ {\lim_{x\to y}}^\kappa f(x), & y\in Y\bs \kappa(X).\end{dcases}\]
    Let $y\in Y$ and take a sequence $( y_n)_{n \in \mathbb{N}} \subset Y$, $y_n \to y$.
        Since $\kappa(X)$ is dense in $Y$, we have that for every $n\in\bN$ there exists $ (x_{n,\,j})_{j \in \mathbb{N}} \subset X$ such that $ \lim_{j\to\infty}\kappa( x_{n,\,j}) = y_n$.
        Since $(\kappa(x_{n,j}))_{j\in\bN}$ converges to $y_n$, for every $n \in \mathbb{N}$, there exists $a_n \in \mathbb{N}$ such that, for every $ j \geq a_n$,
        $d( \kappa(x_{n,\,j}) , y_n) \leq \frac{1}{n}$.

        Now, let us consider two different cases:
        \begin{enumerate}
        	\item If $y_n\in \kappa(X)$, since $\kappa^{-1}|_{\kappa(X)}$ is continuous, we have that $x_{n,j}\to \kappa^{-1}(y_n)$ and, therefore, since $f$ is continuous, there exists $b_n\in\bN$ such that, for $j\ge b_n$,
         \[d( f(x_{n,\,j}) , \widetilde{f}(y_n) )=d(f(x_{n,\,j}) , {f}(\kappa^{-1}(y_n)) ) \leq \frac{1}{n}.\]
         \item On the other hand, if $y_n\in Y\bs\kappa(X)$, then $\til f(y_n):=\lim_{x\to y_n}^\kappa f(x)$ and so  there exists $\d\in\bR^+$ such that $d(f(x),\til f(y_n))<\frac{1}{n}$ if $d(\kappa(x),y_n)<\d$. We have that  $(\kappa(x_{n,j}))_{j\in\bN}$ converges to $y_n$, so there exists $b_n\in\bN$ such that, for $j\ge b_n$, $d(\kappa(x_{n,j}),y_n)<\d$  and, therefore,
        \[d(f(x_{n,\,j}),\widetilde{f}(y_n) ) \leq \frac{1}{n}.\]
    \end{enumerate}
        Hence, for $j\ge j_n:=\max\{a_n,b_n\}$,\[d( \kappa(x_{n,\,j}) , y_n ) \leq \frac{1}{n} \quad \text{ and }\quad  d(f(x_{n,\,j}) , \widetilde{f}(y_n) ) \leq \frac{1}{n}.\]

    Let us define $z_n := x_{n,\,j_n}$ for every $n \in \mathbb{N}$. By the triangle inequality, we have that
    \[d(\kappa( z_n) , y) \leq d(\kappa( z_n)  ,y_n ) + d( y_n , y) \leq \frac{1}{n} + d( y_n , y).\]
    Since $y_n \to y$, we have that $d( y_n, y )\to 0$. Thus,  $d(\kappa( z_n) , y) \to 0$. Then, we are in one of the following cases:
    \begin{enumerate}
    	\item If $y\in Y\bs \kappa(X)$,  $\til f(y)=\lim_{x\to y}^\kappa f(x)$. Then, since $d(\kappa(z_n),y) \rightarrow 0$, $d(f(z_n),\til f(y)) \rightarrow 0$ and, as a consequence, $\lim_{n\to\infty}f(z_n)=\til f(y)$.
    	\item If $y\in \kappa(X)$ then, by the continuity of $\kappa^{-1}|_{\kappa(X)}$, we have that $z_n\to\kappa^{-1}(y)$ and hence, by the continuity of $f$, $\lim\limits_{n \to \infty} f(z_n)=f(\kappa^{-1}(y))=\til f(y)$.
   \end{enumerate}
   In any case, \[\lim\limits_{n \to \infty} f(z_n) = \widetilde f(y).\] Consequently,
    \[d( \widetilde f(y_n) ,\widetilde f(y)) \leq d (\widetilde f(y_n) , f (z_n) ) + d( f(z_n) , \widetilde f(y) ) \leq \frac{1}{n} + d( f(z_n) , \widetilde f(y) ) \to 0.\]
    Therefore, $\lim\limits_{n \to \infty}\widetilde f(y_n) = \widetilde f(y)$.

(II)$\Ra$(I) Let $y\in Y\bs\kappa(X)$ and $V$ a neighborhood of $\til f(y)\in Z$. Let $U=\til f^{-1}(V)$. Since $\til f$ is continuous, $U$ is a neighborhood of $y$ and, since $f=\til f\circ\k$,
\[f(\k^{-1}(U))=f(\k^{-1}(\til f^{-1}(V)))=f(f^{-1}(V))\ss V,\]
so $\lim_{x\to y}^\kappa f(x)=\til f(y)$.

The uniqueness of the extension is due to Proposition \ref{proun}.
\end{proof}
 \section{The space of continuously $\boldmath{n}$-differentiable\\$\boldmath{ (\kappa,\varphi)}$-extensions}

 It is now our objective to study the analytic properties of functions in metric compactifications of closures of open sets in $\bR^n$. The key to achieve this is to contruct an adequate Banach space that we will call $\cC^m_{\kappa,\phi}$.

First, let us introduce some notation. If $\alpha=(\alpha_1,\dots,\alpha_n)\in(\{0\}\cup{\mathbb N})^n$, we define
  \[
  |\alpha|:= \sum_{j=1}^n \alpha_j;
  \]
  and
  \[
  \partial_{\alpha} := \dfrac{\partial^{|\alpha|}}{\partial x_1^{\alpha_1} \cdots \partial x_n^{\alpha_n}}.
  \]
  In particular, for $k \in \{1,\dots,n\}$ and $e_k = (\delta_{1,k},\dots,\delta_{n,k})$, we simply denote
  $\partial_{e_k}$ by $\partial_k$. Let $P_m:=\{p\in(\{0\}\cup{\mathbb N})^n\ :\ |p|\le m\}$.

Let $n,m\in\bN$, $A\subset\bR^n$ open and connected and unbounded, $X$ a compact (and thus complete and totally bounded --see \cite[Theorem 45.1]{Mun75}) metric space, $\kappa:\ol A\to X$ a compactification and $\phi\in\cC^m(\ol A,\bR^+)$. We denote by $\cC(X,\bR)$ the space of continuous functions from $X$ to $\bR$. $\cC(X,\bR)$  is a Banach space with the usual supremum norm: $\|f\|_\infty=\sup_{x\in X}\|f(x)\|_\infty$.

 Define
\[\cC^m_{\kappa,\phi}(\ol A):=\left\{f\in\cC^m(\ol A,\bR):\  \exists {\lim_{y\to x}}^\kappa\partial_p(f/\varphi)(y) \in \bR,\ x\in X\bs\kappa(\ol A),\ p\in P_m\right\}.\]

\begin{lem}\label{lemcar}
\[\cC^m_{\kappa,\phi}(\ol A)=\left\{f\in\cC^m(\ol A,\bR):\  \exists \til f_p\in\cC(X,\bR),\ \partial_p(f/\varphi)=\til f_p\circ\k,\ p\in P_m\right\}.\]
Furthermore, the $\til f_p$ are unique.
\end{lem}
\begin{proof}
 Let $f\in \cC^m_{\kappa,\phi}(\ol A)$, $p\in P_m$. Since for every $x\in X\bs\kappa(\ol A)$ and every $p\in P_m$ there exists $\lim_{y\to x}^\kappa\pd_p(f/\varphi)(x)$, by Theorem~\ref{lemex}, there exists a continuous map $\til f_p\in\cC(X,\bR)$ such that $\pd_p(f/\phi)=\til f_p\circ\k$, so
 \[\cC^m_{\kappa,\phi}(\ol A)\ss\left\{f\in\cC^m(\ol A,\bR):\  \exists \til f_p\in\cC(X,\bR),\ \partial_p(f/\varphi)=\til f_p\circ\k,\ p\in P_m\right\}.\]

 On the other hand, if $f\in\cC^m(\ol A,\bR)$ is such that for every $p\in P_m$ there exists $\til f_p\in\cC(X,\bR)$ satisfying $\partial_p(f/\varphi)=\til f_p\circ\k$, by Theorem~\ref{lemex}, there exists $\lim_{y\to x}^\kappa \partial_p (f/\varphi)(x)$ for every $x\in X\bs\kappa(\ol A)$ and $p\in P_m$.

 Finally, if $g,h\in\cC(X,\bR)$ are such that $g\circ\k=h\circ\k$, since $g$ and $h$ are continuous and $\k(\ol A)$ is dense in $X$, $g=h$. Therefore, the $\til f_p$ are unique.
\end{proof}
\begin{rem} In \cite{Cabada2018} the authors identify the spaces
	\[\left\{f:\bR\to\bR\ :\ f|_{\bR}\in\cC^m(\bR,\bR),\ \exists \lim_{t\to\pm\infty}f^{(j)}(t) \in \bR,\ j=0,\dots,m\right\},\]
	and
	\[\widetilde\cC^m_\varphi:=\left\{f\in\cC^m(\bR,\bR)\ :\ \exists \til f\in\cC^m(\ol\bR,\bR),\  f=\varphi \cdot\til f|_{\bR}\right\},\]
	where $\ol\bR=[-\infty,\infty]$ is the extended real line with its usual topology, observing in \cite[Remark~3.3]{Cabada2018} that if $f\in\cC(\ol \bR,\bR)$ and $f|_\bR\in\cC^m(\bR,\bR)$ then $\lil_{t\to\pm\infty}f^{(j)}(t)=0$ for every $j=1,\dots,m$ since $f$ is asymptotically constant. This is not true in general, as the following example shows. Nonetheless, it is enough to define \[\widetilde\cC^m_\varphi:=\left\{f\in\cC^m(\bR,\bR)\ :\ \exists \til f_j\in\cC^m(\ol\bR,\bR),\  (f/\phi)^{(j)}=\til f_j|_{\bR},\ j=0,\dots,m\right\}.\]
	to obtain the identity.
	\end{rem}
\begin{exa}
Let $g:\bR\to\bR$ be such that $g(x):=(1-x^2)^2$ if $x\in[-1,1]$ and $g(x)=0$ otherwise. Let \[g_k(x):=\frac{g\left(k x-k^2\right)}{k}\] for every $x\in\bR$, $k\ge 2$. Observe that $\supp g_k=[k-\frac{1}{k},k+\frac{1}{k}]$, so $\supp g_k\cap \supp g_j=\emptyset$ for every $k,j\ge 2$, $k\ne j$. Therefore, the function $f(x):=\sum_{k=2}^\infty g_k(x)$, $x\in\bR$, is well defined. The function $f$ and its derivative are illustrated in Figure~\ref{fig:f1}.
\begin{figure}[h]
	\centering
	\includegraphics[width=0.45\linewidth]{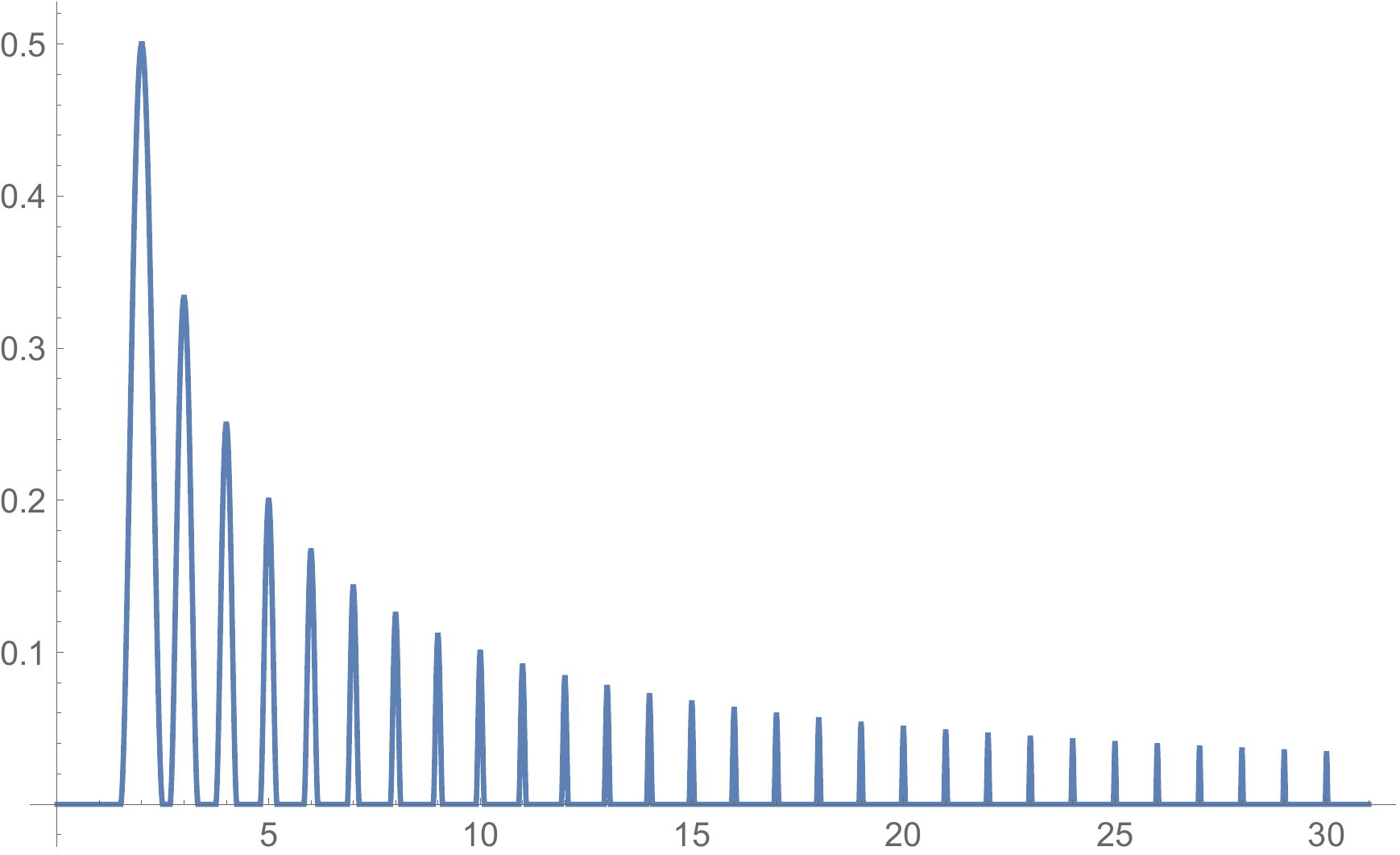}
	\includegraphics[width=0.45\linewidth]{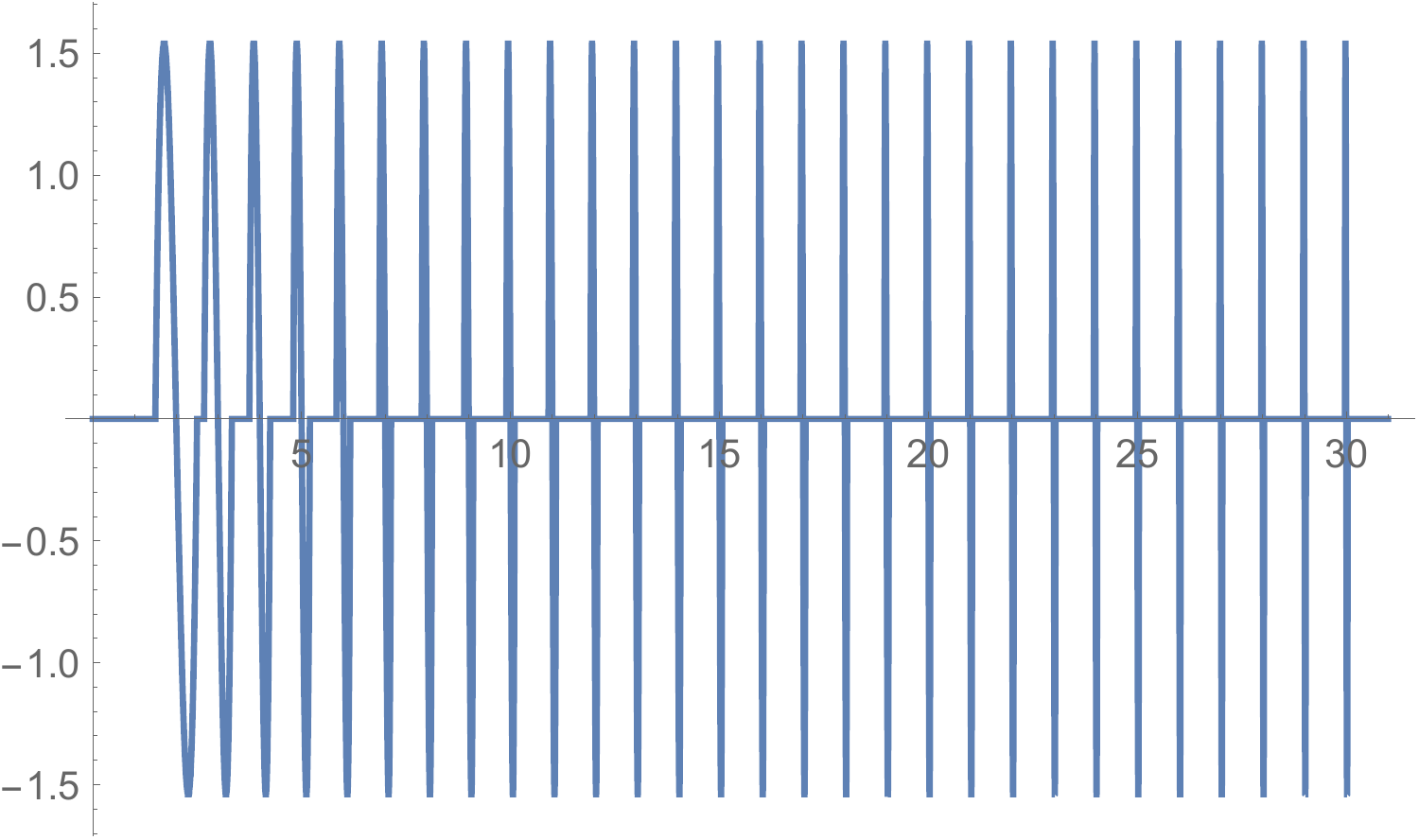}
	\caption{Representation of the function $f$ (left) and its derivative (right).}
	\label{fig:f1}
\end{figure}

Observe that $|g_k(x)|\le\frac{1}{k}$, so $\lim_{x\to\infty}f(x)=0$. Furthermore, since $g_k\in\cC^1(\bR,\bR)$ for every $k\ge 2$, $f\in\cC^1(\bR,\bR)$, but it is not true that there exists $\lim_{x\to\infty}f'(x)$. Indeed: \[f'\(k-\frac{1}{\sqrt{3} k}\)=g_k'\(k-\frac{1}{\sqrt{3} k}\)=\frac{8}{3 \sqrt{3}},\ k\ge 2,\]
but
\[f'\(k+\frac{1}{k}\)=0,\ k\ge 2,\]
so $\lim_{x\to\infty}f'(x)$ does not exist.

\end{exa}

In the next result we will prove that $\cC^m_{\kappa,\phi}(\ol A)$ is a Banach space. To do that we first consider the Banach space $\cB\cC^m(\ol A)$ of $m$-times continuously differentiable  bounded real functions $f$ with the norm \[\|f\|_{m}:=\max\left\{\left\|\partial_p f\right\|_\infty\ :\ p\in P_m\right\}.\]
\begin{thm}\label{t:bound} $\cC^m_{\kappa,\phi}(\ol A)$ is a Banach space with the norm
\[\|f\|_{\kappa,\phi}:=\|f/\varphi\|_{m}, \ f\in \cC^m_{\kappa,\phi}(\ol A).\]
In fact, $\cC^m_{\kappa,\phi}(\ol A)$ is isometrically isomorphic to a closed subspace of $\cB\cC^m(\ol A)$.
\end{thm}
\begin{proof}
Given the linearity of the limits it is clear that $\cC^m_{\kappa,\phi}(\ol A)$ is a vector space.  Consider the map
	\begin{center}\begin{tikzcd}[row sep=tiny]
	\cC^m_{\kappa,\phi}(\ol A)	\arrow{r}{\Xi} & \cB\cC^m(\ol A)  \\
	f\arrow[hook]{r} & f/\phi.
	\end{tikzcd}
\end{center}
\emph{$\Xi$ is well defined:} If $f\in\cC^m_{\kappa,\phi}(\ol A)$, then, for every $p\in P_m$  $\til f_p$ satisfies $\partial_p(f/\varphi)=\til f_p\circ\k$. This means that $f/\varphi$ admits a continuous $p$ derivative for every $p\in P_m$. Since $X$ is compact, $\til f_p$ is bounded, so $\partial_p(f/\varphi)$ is bounded for every $p\in P_m$. Thus, $f/\varphi\in \cB\cC^m(\ol A)$.


$\Xi$ is clearly linear and injective and we can induce the norm of $\cB\cC^m(\ol A)$ in $\cC^m_{\kappa,\phi}(\ol A)$ as $\|f\|_{\kappa,\phi}:=\|\Xi f\|_{m}=\|f/\varphi\|_{m}$. With that norm, $\Xi$ is a bounded isometric map and $\cC^m_{\kappa,\phi}(\ol A)$ and $\Xi(\cC^m_{\kappa,\phi}(\ol A))$ are isometrically isomorphic. It is left to check that $\Xi(\cC^m_{\kappa,\phi}(\ol A))$ is a closed subset of $\cB\cC^m(\ol A)$.

Take a sequence $(f_k)_{k\in\bN}\ss \cC^m_{\kappa,\phi}(\ol A)$ such that $\lim_{k\to 0}\|f_k/\phi- g\|_m=0$ for some $g\in\cB\cC^m(\ol A)$. Let $f=g\phi\in\cC^m(\ol A,\bR)$. Fix $p\in P_m$ and let us check that there exists $h_p\in\cC(X,\bR)$ such that $\partial_pg=\partial_p(f/\varphi)=h_{p}\circ\k$. For every $k\in\bN$ there exists $\til f_{k,p}\in\cC(X,\bR)$ such that $\partial_p(f_k/\varphi)=\til f_{k,p}\circ\k$. We know that $\lim_{k\to 0}\|\partial_p(f_k/\phi)-\partial_p g\|_\infty=0$, so
$\lim_{k\to 0}\|\til f_{k,p}\circ\k-\partial_p g\|_\infty=0$.

Let $\e\in\bR^+$ be fixed. There exists $N\in\bN$ such that $\|\til f_{k,p}\circ\k-\partial_p g\|_\infty<\frac{\e}{2}$ for $k\ge N$. Thus, for $k,j\ge N$, since $\til f_{k,p}$ and $\til f_{j,p}$ are continuous,
\begin{align*}\|\til f_{k,p}-\til f_{j,p}\|_\infty= & \sup_{x\in X}|\til f_{k,p}(x)-\til f_{j,p}(x)|=\sup_{x\in \kappa(\ol A)}|\til f_{k,p}(x)-\til f_{j,p}(x)|=\sup_{y\in \ol A}|\til f_{k,p}(\k(y))-\til f_{j,p}(\k(y))| \\= &\|\til f_{k,p}\circ\k-\til f_{j,p}\circ\k\|_\infty\le \|\til f_{k,p}\circ\k-\partial_p g\|_\infty+\|\partial_p g-\til f_{j,p}\circ\k\|_\infty<\e.\end{align*}

This means $(\til f_{k,p})_{k\in\bN}$ is a Cauchy sequence and, since $\cC(X,\bR)$ is a Banach space, it is convergent to some function $h_p\in\cC(X,\bR)$. Let $M\ge N$ be such that $\|h_p-\til f_{k,p}\|_\infty<\frac{\e}{2}$ for every $k\ge M$. Hence,
\[\|h_p\circ \k-\partial_p g\|_\infty\le\|h_p\circ \k-\til f_{k,p}\circ \k\|_\infty+\|\til f_{k,p}\circ \k -\partial_p g\|_\infty\le \|h_p-\til f_{k,p}\|_\infty+\|\til f_{k,p}\circ \k-\partial_p g\|_\infty<\e.\]
Since $\e$ was fixed arbitrarily, $h_p\circ \k=\partial_p g$. Thus, $f\in\cC^m_{\kappa,\phi}(\ol A)$.


%
%

\end{proof}
%

    Lemma~\ref{lemcar} allows us to define, for every $p\in P_m$, a function
    \begin{equation*}\begin{split}
    \Gamma_p:\cC^m_{\kappa,\phi}(\ol A)&\longrightarrow\cC(X,\bR)	\\
    f & \longmapsto \Gamma_pf,
    \end{split}\end{equation*}
     where $\Gamma_pf$ is the unique function satisfying the equality $\Gamma_pf\circ\k=\pd_p(f/\phi)$.
    \begin{lem}
        $\Gamma_p:\cC^m_{\kappa,\phi}(\ol A)\to\cC(X,\bR)$ is a continuous linear map. Furthermore, $\Gamma_0$ is injective.
    \end{lem}
    \begin{proof}
        $\Gamma_p$  is linear by the linearity of $\partial_p$. To see that it is continuous observe that, due to the density of $\k(\ol A)$ in $X$ and the continuity of the functions $\Gamma_pf$,
        \begin{align*}\|f\|_{\kappa,\phi}= & \max\left\{\left\|\partial_p(f/\phi)\right\|_\infty\ :\ p\in P_m\right\}=\max\left\{\|(\Gamma_pf)\circ\kappa\|_\infty\ :\ p\in P_m\right\}\\ = & \max\left\{\|\Gamma_pf\|_\infty\ :\ p\in P_m\right\}.\end{align*}
        Thus, $\|\Gamma_pf\|_\infty\le\|f\|_{\kappa,\phi}$ for every $f\in \cC^m_{\kappa,\phi}(\ol A)$ and $\Gamma_p$ is continuous.

        $\Gamma_0$ is injective by the uniqueness of the functions $\til f_p$.
    \end{proof}
\begin{rem}
	There is an interesting relation between the norm in $\cC^m_{\kappa,\phi}(\ol A)$ and that of the $\C_p$. Remember that, due to the density of $\k(\ol A)$ in $X$ and the continuity of the functions $\Gamma_pf$,
	\begin{align*}\|f\|_{\kappa,\phi}= \max\left\{\|\Gamma_pf\|_\infty\ :\ p\in P_m\right\}.\end{align*}
	Thus, for any $\e\in\bR^+$ and $f,g\in \cC^m_{\kappa,\phi}(\ol A)$,
	\begin{equation}\label{eqballc}g\in B_{\cC^m_{\kappa,\phi}(\ol A)}(f,\e)\iff \C_p g\in B_{\cC(X,\bR)}(\C_p f,\e)\sfa p\in P_m.
	\end{equation}
\end{rem}

In order to successfully develop the next secction we need a precompactness criterion for subsets in $\cC^m_{\kappa,\phi}(\ol A)$. Unfortunately, Ascoli-Arzelà Theorem, as normally stated, cannot be applied to  $\cC^m_{\kappa,\phi}(\ol A)$ nor $\cB\cC^m(\ol A)$ directly as $\ol A$ is not compact. Nevertheless, the theorem does apply to $\cC(X,\bR)$ since $X$ is a Hausdorff compact topological space and $\bR$ is a complete metric space.
\begin{thm}[Ascoli-Arzel\`a \cite{kelley}]\label{aa} Let $X$ be a compact  Hausdorff  topological space and $Y$ a complete metric space, and consider $\cC(X,Y)$ with the topology of the uniform convergence. Then $F\subset\cC(X,Y)$ has compact closure if and only if for every $x\in X$
\begin{enumerate}
\item \emph{$F$ is uniformly bounded at $x$}, that is $F(x):=\{f(x)\ :\ f\in X\}$ has compact closure, and
\item \emph{$F$ is equicontinuous at $x$}, that is, for every $\e\in\bR^+$ there exists a neighbourhood $U$ of $x$ such that $d(f(y),f(x))<\e$ for every $y\in U$ and $f\in F$.
\end{enumerate}
\end{thm}

The question is, what is the relation between compactness in $\cC^m_{\kappa,\phi}(\ol A)$ and in $\cC(X,\bR)$? The answer to this question comes from Lemma~\ref{lemcar}, as the following result shows.

\begin{lem}\label{lemaa}Let $F\subset \cC^m_{\kappa,\phi}(\ol A)$. $F$ is compact if and only if $\Gamma_p(F)$ is compact in $\cC(X,\bR)$ for every $p\in P_m$.
\end{lem}
\begin{proof}
%
%
Let $p\in P_m$ and assume $F$ is compact. Since $\Gamma_p$ is continuous, $\Gamma_p(F)$ is compact.

On the other hand, assume $\Gamma_p(F)$ is compact in $\cC(X,\bR)$ for every $p\in P_m$. Write $P_m=\{p_j\}_{j=1}^r$. Let $\cU$ be an open cover of $F$. For every $f\in F$ there exists $U_f\in\cU$ such that $f\in U_f$. Since $U_f$ is open, there exists $\d_f\in\bR^+$ such that $B_{\cC^m_{\kappa,\phi}(\ol A)}(f,\d_f)\ss U_f$. Let $\cV_{p_1}=\{B_{\cC(X,\bR)}(\C_{p_1} f,\d_f)\}_{f\in F}$. We have that $\cV_{p_1}$ is an open cover of $\Gamma_{p_1}(F)$, which is compact,  so there exists an open subcover $\{B_{\cC(X,\bR)}(\C_{p_1} f_{1,j},\d_{1,j})\}_{j=1}^{h_1}$ of $\cV_{p_1}$.

 Now, for every $f\in F$, let \[J_{1,f}:=\{j\in\{1,\dots,h_1\}\ :\ \C_{p_1}f\in B_{\cC(X,\bR)}(\C_{p_1} f_{1,j},\d_{1,j})\}.\] Because $\cV_{p_1}$ covers $\C_{p_1}(F)$, we know that $J_{1,f}\ne \emptyset$, so there exists $\d_{f,1}\in\bR^+$ such that \[B_{\cC^m_{\kappa,\phi}(\ol A)}(f,\d_{f,1})\ss U_f\cap\bigcap_{j\in F_{1,f}} B_{\cC^m_{\kappa,\phi}(\ol A)}\( f_{1,j},\frac{\d_{1,j}}{2}\).\]  $\cV_{p_2}=\{B_{\cC(X,\bR)}(\C_{p_2} f,\d_{f,1})\}_{f\in F}$ is an open cover of $\Gamma_{p_2}(F)$, which is compact, so there exists an open subcover $\{B_{\cC(X,\bR)}(\C_{p_2} f_{2,j},\d_{2,j})\}_{j=1}^{h_2}$ of $\cV_{p_2}$.

We repeat this process: having constructed $\{B_{\cC(X,\bR)}(\C_{p_1} f_{k,j},\d_{k,j})\}_{j=1}^{h_{k}}$, an open subcover of $\cV_{p_k}$, we define, for every $f\in F$, \[J_{k,f}:=\{j\in\{1,\dots,h_k\}\ :\ \C_{p_k}f\in B_{\cC(X,\bR)}(\C_{p_k} f_{k,j},\d_{k,j})\}\ne\emptyset,\] so there exists $\d_{f,k}\in\bR^+$  such that \[B_{\cC^m_{\kappa,\phi}(\ol A)}(f,\d_{f,k})\ss U_f\cap\bigcap_{l=1}^k \bigcap_{j\in F_{l,f}} B_{\cC^m_{\kappa,\phi}(\ol A)}\( f_{l,j},\frac{\d_{l,j}}{2}\).\]
Then $\cV_{p_{k+1}}=\{B_{\cC(X,\bR)}(\C_{p_k} f,\d_{f,k})\}_{f\in F}$ is an open cover of $\Gamma_{p_{k+1}}(F)$, which is compact, so there exists an open subcover $\{B_{\cC(X,\bR)}(\C_{p_{k+1}} f_{k+1,j},\d_{k+1,j})\}_{j=1}^{h_{k+1}}$ of $\cV_{p_{k+1}}$.

$\{B_{\cC(X,\bR)}(\C_{p_r} f_{r,j},\d_{r,j})\}_{j=1}^{h_{r}}$ is an open cover of $\cV_{p_{r}}$ and, for every $s=1,\dots,h_{r}$,  \begin{equation}\label{eqball}B_{\cC^m_{\kappa,\phi}(\ol A)}( f_{r,s},\d_{r,s})\ss U_{ f_{r,s}}\cap\bigcap_{l=1}^{r-1} \bigcap_{j\in F_{l, f_{r,s}}} B_{\cC^m_{\kappa,\phi}(\ol A)}\( f_{l,j},\frac{\d_{l,j}}{2}\).\end{equation}
Consider the family $\til\cU:=\{U_{f_{r,s}}\}_{j=1}^{h_{r}}$. $\til\cU$ is a finite subset of $\cU$. Let us check that it is a subcover. Take $g\in F$. Then $\C_{p_r}g\in B_{\cC(X,\bR)}(\C_{p_r} f_{r,s},\d_{r,s}/2)$ for some $s=1,\dots, h_r$, so $\|\C_{p_r}g-\C_{p_r} f_{r,s}\|_\infty<\d_{r,s}$. We now want to prove that $\|\C_qg-\C_qf_{r,s}\|_\infty<\d_{r,s}$ for every $q\in\{1,\dots, r-1\}$ because then, by expression~\eqref{eqballc}, $\|g-f_{r,s}\|_{\k,\phi}<\d_{r,s}<$ and, thus, $g\in B_{\cC^m_{\kappa,\phi}(\ol A)}(f_{r,s},\d_{r,s})\ss U_{f_{r,s}}\in \til\cU$.

By expression~\eqref{eqball}, $\|g-f_{l,j}\|_{\k,\phi}<\d_{l,j}/2$ and $\|f_{r,s}-f_{l,j}\|_{\k,\phi}<\d_{l,j}/2$ for every $j\in F_{l, f_{r,s}}$ and $l\in\{1,\dots,r-1\}$, which, by expression~\eqref{eqball}, implies \[\|\C_qg-\C_qf_{l,j}\|_{\infty},\ \|\C_qf_{l,j}-\C_qf_{r,s}\|_{\infty}<\frac{\d_{l,j}}{2},\] for every $q\in P_m$, $j\in F_{l, f_{r,s}}$ and $l\in\{1,\dots,r-1\}$. Hence,
\[\|\C_qg-\C_qf_{r,s}\|_{\infty}\le\|\C_qg-\C_qf_{l,j}\|_{\infty}+\|\C_qf_{l,j}-\C_qf_{r,s}\|_{\infty}<\d_{l,j}<\d_{r,s},\]
as we wanted to show.
%

\end{proof}

Now using Lemma~\ref{lemaa}, it is possible to rewrite Theorem \ref{aa} in terms of $\cC^m_{\kappa,\phi}(\ol A)$, getting the following result.

\begin{thm}\label{taae}
Let $A\subset\bR^n$ be open and unbounded, $X$ be a compact metric space and $\kappa:\ol A\to X$ a compactification. $F\subset \cC^m_{\kappa,\phi}(\ol A)$ has compact closure if the three following conditions are satisfied:
\begin{enumerate}
	\item  For each $x\in\ol A$ and every $p\in P_m$ there exists some constant $M_{x,p}\in\bR^+$ such that
	\[\left|\partial_p(f/\varphi)(x)\right|\le M_{x,p},\]
	for all $f\in F$.
	\item For every $x\in \ol A$, $\varepsilon\in \bR^+$ and $p\in P_m$ there exists some $\delta_{x,p}\in\bR^+$ such that
	\[\left|\partial_p(f/\varphi)(x)-\partial_p(f/\varphi)(y)\right|<\varepsilon,\]
	if $f\in F$ and $y\in\ol A$ is such that $\|x-y\|<\delta_{x,p}$.
    \item For every $x\in X\bs\k(\ol A)$, $\varepsilon\in \bR^+$ and $p\in P_m$ there exists some $\delta_{x,p}\in\bR^+$ such that
    \[\left|{\lim_{z\to x}}^\kappa\partial_p(f/\varphi)(z)-\partial_p(f/\varphi)(y)\right|<\varepsilon,\]
    if $f\in F$ and $y\in\ol A$ is such that $d(x,\k(y))<\delta_{x,p}$.
\end{enumerate}
\end{thm}
\begin{proof}
Assume 1, 2 and 3 hold.  Fix $p\in P_m$ and $x\in X$.

\textbf{Step 1:} We will show that $\Gamma_p(F)$ is equicontinuous at $x$. We study two cases:
\begin{enumerate}
	\item[\emph{(a)}]
 If $x\in\k(\ol A)$, then $x=\k(z)$ for some $z\in\ol A$ and for every $\varepsilon \in \bR^+$ there exists $\d_{x,p}$ such that for every $f\in F$, $y\in\ol A$, $\|y-z\|<\d_{x,p}$,
\[\left|\Gamma_pf(x)-\Gamma_pf(\k(y))\right|=\left|\Gamma_pf(\k(z))-\Gamma_pf(\k(y))\right|=\left|\partial_p(f/\varphi)(z)-\partial_p(f/\varphi)(y)\right|<\frac{\varepsilon}{2}.\]
Since $\k:\ol A\to\k(\ol A)$ is an homeomorphism, $\k^{-1}:\k(\ol A)\to\ol A$ is continuous, so there  exists $\til\d_{x,p}\in\bR^+$ such that if $\k(y)\in \k(\ol A)$ and $d(\k(y),x)<\d_{x,p}$ then $\|z-y\|<\til\delta_{x,p}$.
Thus
\[\Gamma_pf\left(B_X(x,\til\delta_{x,p})\cap\k(\ol A)\right)\ss B_\bR\(\Gamma_pf(x),\frac{\varepsilon}{2}\).\]
Since $\k(\ol A)$ is dense in $X$ and $\Gamma_p$ is continuous,
\[\Gamma_pf\left(B_{X}(x,\til\delta_{x,p})\right)\ss B_\bR\left[\Gamma_pf(x),\frac{\varepsilon}{2}\right]\ss B_\bR\left(\Gamma_pf(x),\e\right)\]
for every $f\in F$, which implies that $\Gamma_p(F)$ is equicontinuous at $x$ in $\cC(X,\bR)$.

\item[\emph{(b)}]  If $x\in X\bs\k(\ol A)$, we resort to an analogous argument using 3 instead of 2. There exists $\d_{x,p}$ such that for every $f\in F$, $y\in X$, $d(x,\k(y))<\delta_{x,p}$,
\[\left|\Gamma_pf(x)-\Gamma_pf(\k(y))\right|=\left|{\lim_{z\to x}}^\kappa\partial_p(f/\varphi)(z)-\partial_p(f/\varphi)(y)\right|<\frac{\varepsilon}{2},\]
by the continuity of $\Gamma_p$. Since $\k:\ol A\to\k(\ol A)$ is a homeomorphism, $\k^{-1}:\k(\ol A)\to\ol A$ is continuous, so there exists $\til\d_{x,p}\in\bR^+$ such that if $\k(y)\in \k(\ol A)$ and $d(\k(y),x)<\d_{x,p}$ then $\|z-y\|<\til\delta_{x,p}$.
Thus
\[\Gamma_pf\left(B_X(x,\til\delta_{x,p})\cap\k(\ol A)\right)\ss B_\bR\(\Gamma_pf(x),\frac{\varepsilon}{2}\).\]
Since $\k(\ol A)$ is dense in $X$ and $\Gamma_p$ is continuous,
\[\Gamma_pf\left(B_{X}(x,\til\delta_{x,p})\right)\ss B_\bR\left[\Gamma_pf(x),\frac{\varepsilon}{2}\right]\ss B_\bR\left(\Gamma_pf(x),\e\right)\]
for every $f\in F$, which implies that $\Gamma_p(F)$ is equicontinuous at $x$ in $\cC(X,\bR)$.
\end{enumerate}

\textbf{Step 2:} We  will show that $\Gamma_p(F)$ is uniformly bounded at $x$. We again have two cases:
\begin{enumerate}
\item[\emph{(a)}]  If $x=\kappa(z)$ for some $z\in \ol A$, by 1, we have that $|\C_pf(x)|=|\C_pf(\k(z))|\le M_{x,p}$ for all $f\in F$ and that is it.
\item[\emph{(b)}]  Assume now $x\in X\bs\k(\ol A)$.
$\C_p(F)$ is equicontinuous at $x$ so  there exists $\d_{x,p}\in\bR^+$ such that, for every $f\in F$,
\[\Gamma_pf\left(B_{X}(y,\delta_{x,p})\right)\ss B_\bR\left(\Gamma_pf(x),1\right).\]
Take $\k(t)\in B_{X}(y,\delta_{x,p})\cap\k(\ol A)$.
By 1, we have that $|\C_pf(\k(t))|\le M_{t,p}$ for every $f\in F$. Thus, for every $f\in F$,
\[|\Gamma_pf(x)|\le|\Gamma_pf(x)-\C_pf(\k(t))|+|\C_pf(\k(t))|\le 1+ M_{t,p}.\]
\end{enumerate}
Hence, $\Gamma_pF$ is uniformly bounded at $x$.

We conclude that $\ol{\Gamma_pF}$ is compact for every $p\in P_m$ and, thus, $\ol F$ is compact.

Finally, assume $\ol F$ is compact. Then so is $\ol{\Gamma_pF}$ for every $p\in P_m$. Fix $p\in P_m$ and $x\in X$. $\Gamma_pF$ is uniformly bounded at $x$, so there exist $ M_{x,p}\in\bR^+$ such that, for every $f\in F$,
\[\left|\C_pf(x)\right|\le M_{x,p}.\]
Thus, given $x\in\ol A$, for every $f\in F$,
\[\left|\partial_p(f/\varphi)(x)\right|=\left|\C_pf(\k(x))\right|\le M_{\k(x),p},\]
so 1 holds.

$\Gamma_pF$ is also equicontinuous at $x$. So for every $\e\in\bR^+$ there exists $\d_{x,p}\in\bR^+$ such that, for every $f\in F$,
\[\Gamma_pf(B_{X}(x,\delta_{x,p}))\ss B_\bR\left(\Gamma_pf(x),\e\right).\]
If $x=\k(z)$ with $z\in \ol A$, then, for every $f\in F$, $y\in\ol A$, $\|y-z\|<\d_{x,p}$,
\[\left|\Gamma_pf(x)-\Gamma_pf(\k(y))\right|=\left|\Gamma_pf(\k(z))-\Gamma_pf(\k(y))\right|=\left|\partial_p(f/\varphi)(z)-\partial_p(f/\varphi)(y)\right|<\e,\]
and 2 holds.
If $x\in X\bs\k(\ol A)$, then, by the continuity of $\C_p$,  for every $f\in F$, $y\in X$, $d(x,\k(y))< \delta_{x,p}$,
\[\left|\Gamma_pf(x)-\Gamma_pf(\k(y))\right|=\left|{\lim_{z\to x}}^\kappa\partial_p(f/\varphi)(z)-\partial_p(f/\varphi)(y)\right|<\e,\]
and 3 holds.
\end{proof}
\begin{rem}Conditions 2 and 3 in Theorem~\ref{taae} can be synthesised as a single one:
 \begin{enumerate}
     \item[4.] For every $x\in X$, $\varepsilon\in \bR^+$ and $p\in P_m$ there exist some $\delta_{x,p}\in\bR^+$ such that
     \[\left|{\lim_{z\to x}}^\kappa\partial_p(f/\varphi)(z)-\partial_p(f/\varphi)(y)\right|<\varepsilon,\]
     if $f\in F$ and $y\in\ol A$ is such that $d(x,\k(y))<\delta_{x,p}$.
 \end{enumerate}
In practice it is convenient to keep them separated as 2 avoids the direct use of the compactification.
\end{rem}
\begin{rem} Condition 3 in Theorem~\ref{taae} is necessary, as the following example shows.
   \end{rem}
\begin{exa}\label{exadest} Consider the extended real line $\ol\bR:=[-\infty,\infty]$ with the topology given by the distance
    \begin{align*}d(x,y)= & \left|\frac{x}{1+x^2}-\frac{y}{1+y^2}\right|,\quad x,y\in\bR,\\ d(\infty,x)= & d(x,\infty)=\left|\frac{x}{1+x^2}-1\right|,\quad d(-\infty,x)=d(x,-\infty)=\left|\frac{x}{1+x^2}+1\right|,\quad x\in\bR.\end{align*}
     $(\ol\bR,d)$ is a compact metric space.  Consider the compactification $\k:\bR\to\ol\bR$ defined as $\k(x)=x$ for every $x\in\bR$. Take $\phi(x)=x$ for every for every $x\in\bR$ and consider the family $F:=\{f_n\}_\n$ where $f_n(x)=e^{-(x-n)^2}$, $x\in\bR$.  $F\ss\til\cC^1_{\kappa,\phi}(\bR)$ is equicontinuous and uniformly bounded (in $\bR$), so it satisfies conditions 1 and 2 in Theorem~\ref{taae}, but not 3, as $\Gamma_0(F)\ss\cC(\ol\bR,\bR)$ is not compact. Indeed: $\|f_n-f_m\|\ge 1-e^{-1}$ for every $m,n\in\bN$, $m\ne n$. This means that $F$ admits no Cauchy sequence and, thus, no convergent subsequence, so $F$ cannot be compact.
\end{exa}
\begin{rem}
We observe that in \cite[Theorem 3.2]{Cabada2018} condition 3 of Theorem~\ref{taae} is missing. Also, \cite[Theorem 1]{przerad} can be considered as a particular instance of Theorem~\ref{taae} for the case of the compactification in Example~\ref{exadest}. Condition 3 of Theorem~\ref{taae} is called \emph{regularity} in \cite{przerad}. A similar condition appears in \cite[Section 2.12, p. 62]{Corduneanu} under the name \emph{equiconvergence} in a setting that would correspond to the compactification $\til\k=\k|_{[0,\infty)}$ where $\k$ is taken as in Example~\ref{exadest}.
    \end{rem}
%

%

\section{Fixed points of integral equations}
In this section we will prove the existence of fixed points of integral equations. To that end, we will develop a method based on the fixed point index theory on abstract cones.

With the notation introduced in the previous section, let $n,\,m\in\bN$, $A\subset \bR^n$ open, connected and unbounded, $X$ a compact metric space, $\k \colon \ol A \to X$ a compactification and $\varphi \in \cC^m(\ol A, \bR^+)$. Let us consider the integral operator $T\colon \cC_{\k,\varphi}^m (\ol A) \to \cC_{\k,\varphi}^m (\ol A)$, given by the following expression
\[Tu(t)=\int_{\ol A} G(t,s)\,f(s,u(s))\, \dif s, \]
for $t\in \ol A$ (it is defined as the limit of such expressions on $X\bs\k(\ol A)$),
where $G\colon \ol A\times \ol A \rightarrow \bR$. Since the kernel $G$ is defined on $\ol A\times \ol A$, we need to introduce the following notation: given $p\in \left(\{0\}\cup \bN\right)^n$, and using the natural injection
\begin{equation*}\begin{split}
\left(\{0\}\cup \bN\right)^n & \longrightarrow \left(\{0\}\cup \bN\right)^{2n} \\
p=(p_1,\dots,p_n) & \longmapsto \widetilde{p}=(p_1,\dots,p_n,0,\dots,0),
\end{split}\end{equation*}
we shall denote $\partial_{\widetilde{p}} G$ by $\partial_p G(t,\cdot)$. This notation is not to be confused with $\partial_p(G(t,\cdot))$, where we fix the value of $t$ and differentiate with respect to the rest of the derivatives. The same applies to $\partial_p(G(\cdot,s))$ for a fixed $s$.

The next result will provide sufficient conditions for the operator  $T\colon \cC_{\k,\varphi}^m (\ol A) \to \cC_{\k,\varphi}^m (\ol A)$ to be well defined, continuous and compact. We will make use of the following hypotheses:
\begin{itemize}
\item[$(C_1)$] The kernel $G\colon \ol A\times \ol A \rightarrow \bR$ is such that $G(\cdot,s)\in \cC_{\k,\varphi}^m (\ol A)$ for every $s\in \ol A$.

In particular, from the definition of $\cC_{\k,\varphi}^m (\ol A)$ and Theorem~\ref{t:bound}, this implies that there exist
\[M_p(s) :=\sup_{t\in \ol A} \left| \partial_p\left(\frac{G(\cdot,s)}{\varphi}\right)(t) \right|\in\bR, \quad \forall\, s\in\ol A, \  p\in P_m \]
and
\[z_p^x(s):= {\lim\limits_{t\to x}}^{\k} \partial_p\left(\frac{G(\cdot,s)}{\varphi}\right)(t) \in \bR, \quad \forall\, s\in\ol A, \ x\in X\setminus\kappa(\ol A), \  p\in P_m. \]
Taking into account the definition of function $\Gamma_p f$ as the unique function satistying the equality $\Gamma_pf\circ \kappa =\partial_p(f/\varphi)$, and the proof of Theorem~\ref{lemex}, it occurs that \[z_p^x(s)=\Gamma_p\left(G(\cdot,s)\right)(x).\]
\item[$(C_2)$] For every $\varepsilon\in \bR^+$ and $p\in P_m$, there exist $\delta\in \bR^+$ and a measurable function $w_p$ such that, if $x,\,y\in X$ satisfy $d(x,y)<\delta$, then
\begin{equation*}
\left| \Gamma_p\left(G(\cdot,s)\right)(x) - \Gamma_p\left(G(\cdot,s)\right)(y) \right| <	\varepsilon\, w_p(s), \quad \forall\,s\in \ol A.
\end{equation*}
\item[$(C_3)$] The nonlinearity $f\colon \bR^n\times \bR \to [0,\infty)$ satisfies the following conditions:
\begin{itemize}
	\item[\textbullet] $f(\cdot,y)$ is measurable for each fixed $y\in\bR$.
	\item[\textbullet] $f(t,\cdot)$ is continuous for a.\,e. $t\in \bR^n$.
	\item[\textbullet] For each $r>0$ there exists $\Phi_r\in \Lsp{1}(\ol A)$ such that
	\[f(t,y\,\varphi(t))\le \Phi_r(t),\]
	for all $y\in\bR$ with $|y|<r$ and a.\,e. $t\in \bR^n$.
\end{itemize}
\item[$(C_4)$] For every $r>0$, $x\in X\setminus\kappa(\ol A)$ and $p\in P_m$ it holds that $M_p\,\Phi_r, \, |z_p^x|\,\Phi_r, \, w_p\,\Phi_r \in \Lsp{1}(\ol A)$.
\end{itemize}

\begin{thm}
Assume that hypotheses $(C_1)$--$(C_4)$ are satisfied. Then operator $T$ is well defined, continuous and compact.
\end{thm}
\begin{proof}
We shall divide the proof into several steps.

\textbf{Step 1:} Let us prove first that operator $T$ is well defined, that is, that it maps $\cC_{\k,\varphi}^m (\ol A)$ to $\cC_{\k,\varphi}^m (\ol A)$.

From the general rules of differentiability of integrals (see \cite[Corollary~2.8.7]{Bogachev}) it holds that
\begin{equation}\label{eq:part_int}\begin{split}
\partial_p \left(\frac{Tu}{\varphi}\right) (t) &= \partial_p \int_{\ol A} \frac{G(t,s)}{\varphi(t)}\,f(s,u(s))\, \dif s = \int_{\ol A} \partial_p\left(\frac{G(\cdot,s)}{\varphi}\right)(t)\,f(s,u(s))\, \dif s.
\end{split}\end{equation}

Now, fix $\varepsilon\in \bR^+$. From $(C_2)$, there exists $\delta\in \bR^+$ such that if $t_1,\,t_2\in \ol A$, $\|t_1-t_2\|<\delta$ then
\begin{equation*}
	\left| \Gamma_p\left(G(\cdot,s)\right)(\kappa(t_1)) - \Gamma_p\left(G(\cdot,s)\right)(\kappa(t_2)) \right| <	\varepsilon\, w_p(s), \quad \forall\,s\in \ol A,
\end{equation*}
or, which is the same,
\begin{equation*}
\left|	\partial_p \left(\frac{G(\cdot,s)}{\varphi}\right)(t_1) - \partial_p \left(\frac{G(\cdot,s)}{\varphi}\right)(t_2) \right| <\varepsilon\, w_p(s), \quad \forall\,s\in \ol A.
\end{equation*}
Hence,
\begin{equation*}\begin{split}
\left|\partial_p \left(\frac{Tu}{\varphi}\right) (t_1) - \partial_p \left(\frac{Tu}{\varphi}\right) (t_2) \right|	& \le \int_{\ol A} \, \left| 	\partial_p \left(\frac{G(\cdot,s)}{\varphi}\right)(t_1) - \partial_p \left(\frac{G(\cdot,s)}{\varphi}\right)(t_2)	 \right|	\, f(s,u(s))\, \dif s \\
& \le \varepsilon \int_{\ol A} w_p(s)\, f(s,u(s))\, \dif s \le \varepsilon \int_{\ol A} w_p(s)\, \Phi_{\|u\|_{\k,\, \varphi} }(s)\, \dif s.
\end{split}\end{equation*}
Now, from $(C_4)$, it is ensured the existence of some positive constant $c$ such that the previous expression is upperly bounded by $\varepsilon\,c$. Consequently, $\partial_p \left(\frac{Tu}{\varphi}\right)$ is continuous on $\ol A$ for every $p\in P_m$, that is, $\frac{Tu}{\varphi}\in \cC^m(\ol A, \bR)$ and, since $\varphi \in \cC^m(\ol A, \bR^+)$, we conclude that $Tu \in \cC^m(\ol A, \bR)$.

Let us show now that for every $x\in X\setminus\k(\ol A)$ there exists ${\lim\limits_{t\to x}}^{\k} \partial_p\left(\frac{Tu}{\varphi}\right)(t)$. We have that
\begin{equation*}\begin{split}
{\lim\limits_{t\to x}}^{\k} \partial_p\left( \frac{Tu}{\varphi}\right) (t)= {\lim\limits_{t\to x}}^{\k} \int_{\ol A} \partial_p\left(\frac{G(\cdot,s)}{\varphi}\right)(t)\,f(s,u(s))\, \dif s.
\end{split}\end{equation*}
Now, since
\[\left|\partial_p\left(\frac{G(\cdot,s)}{\varphi}\right)(t)\,f(s,u(s)) \right| \le M_p(s)\, \Phi_{\|u\|_{\k,\, \varphi}}(s) \]
and, from $(C_4)$, $M_p\, \Phi_{\|u\|_{\k,\, \varphi}} \in \Lsp{1}(\ol A)$, by Lebesgue's Dominated Convergence Theorem, we obtain that
\begin{equation}\label{eq:lim_part}\begin{aligned}
		{\lim\limits_{t\to x}}^{\k} \partial_p\left( \frac{Tu}{\varphi}\right) (t)&= {\lim\limits_{t\to x}}^{\k} \int_{\ol A} \partial_p\left(\frac{G(\cdot,s)}{\varphi}\right)(t)\,f(s,u(s))\, \dif s \\
		&= \int_{\ol A} \,  {\lim\limits_{t\to x}}^{\k}  \left( \partial_p\left(\frac{G(\cdot,s)}{\varphi}\right)(t)\right)\,f(s,u(s))\, \dif s  \int_{\ol A} \,  z_p^x(s)\,f(s,u(s))\, \dif s.
\end{aligned}\end{equation}
Thus, since
\begin{equation*}
\left|\int_{\ol A} \,  z_p^x(s)\,f(s,u(s))\, \dif s\right|	\le \int_{\ol A} \,  |z_p^x(s)|\,\Phi_{\|u\|_{\k,\, \varphi}}(s) \, \dif s
\end{equation*}
and, from $(C_4)$, $|z_p^x|\, \Phi_{\|u\|_{\k,\, \varphi}} \in \Lsp{1}(\ol A)$, we have proved the existence of ${\lim\limits_{t\to x}}^{\k} \partial_p\left(\frac{Tu}{\varphi}\right)(t)$ for every ${x\in X\setminus \kappa(\ol A)}$.

Therefore, we conclude that $Tu \in \cC_{\k,\varphi}^m (\ol A)$.

\textbf{Step 2:} Continuity:

Let $(u_n)_{n\in\bN} \subset \cC_{\k,\varphi}^m (\ol A)$ be a sequence which converges to $u$ in $\cC_{\k,\varphi}^m (\ol A)$ and let us show that $(Tu_n)_{n\in\bN}$ converges to $Tu$ in $\cC_{\k,\varphi}^m (\ol A)$.

The convergence of $(u_n)_{n\in\bN}$ to $u$ in $\cC_{\k,\varphi}^m (\ol A)$ implies that, in particular, $u_n(s) \to u(s)$ for a.\,e. $s\in\ol A$ and so, from $(C_3)$, $f(s,u_n(s))\to f(s,u(s))$ for a.\,e. $s\in\ol A$.

Following similar arguments to the ones above, we have that, for every $p\in P_m$,
\begin{equation*}\begin{split}
\left|\partial_p \left(\frac{Tu_n}{\varphi}\right) (t) -\partial_p\left(\frac{Tu}{\varphi}\right) (t) \right| &
\le  \int_{\ol A} \left|\partial_p\left(\frac{G(\cdot,s)}{\varphi}\right)(t)\right|\, |f(s,u_n(s))-f(s,u(s))|\, \dif s \\[2pt]
& \le \int_{\ol A} M_p(s)\, |f(s,u_n(s))-f(s,u(s))|\, \dif s.
\end{split}\end{equation*}
Moreover, the convergence of $(u_n)_{n\in\bN}$ ensures the existence of some positive constant $R$ such that $\|u_n\|_{\kappa,\varphi} \le R$ for every $n\in \bN$, which guarantees that the previous integral is upperly bounded by $2\,\int_{\ol A}\, M_p(s)\,\Phi_R(s) \dif s \in \bR$. Therefore, by Lebesgue's Dominated Convergence Theorem (which can be used because of the continuity of $\Gamma_p$), we obtain that
\begin{equation*}\begin{split}
\lim\limits_{n\to \infty} \, \left\| \partial_p \left(\frac{Tu_n}{\varphi}\right) -\partial_p\left(\frac{Tu}{\varphi}\right) \right\| &
\le \lim\limits_{n\to \infty} \int_{\ol A} M_p(s)\, |f(s,u_n(s))-f(s,u(s))|\, \dif s \\
& = \int_{\ol A} \, \lim\limits_{n\to \infty} M_p(s)\, |f(s,u_n(s))-f(s,u(s))|\, \dif s=0.
\end{split}\end{equation*}
This way we have proved that $(Tu_n)_{n\in\bN}$ converges to $Tu$ in $\cC_{\k,\varphi}^m (\ol A)$. Hence, $T$ is a continuous operator.

\textbf{Step 3:} Compactness:

Let us consider a bounded set $B\subset \cC_{\k,\varphi}^m (\ol A)$, that is, such that there exists a positive constant $R$ for which $\|u\|_{\kappa,\varphi}\le R$ for every $u\in B$. We shall prove that $T(B)$ is relatively compact in $\cC_{\k,\varphi}^m (\ol A)$.

Reasoning as above, we obtain that given $p\in P_m$,
\begin{equation*}\begin{split}
\left|\partial_p\left(\frac{Tu}{\varphi}\right) (t) \right| & \le \int_{\ol A} \, \left|\partial_p\left(\frac{G(\cdot,s)}{\varphi}\right)(t)\right|\, |f(s,u(s))|\, \dif s  \le \int_{\ol A} M_p(s)\, \Phi_R(s)\, \dif s
\end{split}\end{equation*}
for every $t\in \ol A$ and $u\in B$. Therefore, we deduce that $T(B)$ is uniformly bounded.

On the other hand, as we have shown in Step 1, for every fixed $\varepsilon \in \bR^+$ there exists $\delta\in \bR^+$ such that if $t_1,\,t_2\in \ol A$, $\|t_1-t_2\|<\delta$ then
\begin{equation*}\begin{split}
		\left|\partial_p \left(\frac{Tu}{\varphi}\right) (t_1) - \partial_p \left(\frac{Tu}{\varphi}\right) (t_2) \right|	& \le \varepsilon \int_{\ol A} w_p(s)\, f(s,u(s))\, \dif s \le \varepsilon \int_{\ol A} w_p(s)\, \Phi_R(s)\, \dif s,
\end{split}\end{equation*}
and, since $w_p\,\Phi_R\in \Lsp{1}(\ol A)$, we can ensure the existence of some constant $c$ such that
\begin{equation*}
		\left|\partial_p \left(\frac{Tu}{\varphi}\right) (t_1) - \partial_p \left(\frac{Tu}{\varphi}\right) (t_2) \right| <\varepsilon\, c \quad \text{for all } u\in B.
\end{equation*}

Finally, consider $x\in X\setminus\kappa(\ol A)$ and $p\in P_m$. From \eqref{eq:part_int} and \eqref{eq:lim_part}, we obtain that
\begin{equation*}\begin{split}
\left|{\lim_{\tau\to x}}^\kappa\partial_p\left(\frac{Tu}{\varphi}\right)(\tau)-\partial_p\left(\frac{Tu}{\varphi}\right)(t)\right| &=\left|\int_{\ol A} \,  z_p^x(s)\,f(s,u(s))\, \dif s - \int_{\ol A} \partial_p\left(\frac{G(\cdot,s)}{\varphi}\right)(t)\,f(s,u(s))\, \dif s \right| \\[2pt]
& \le  \int_{\ol A} \, \left| z_p^x(s)- \partial_p\left(\frac{G(\cdot,s)}{\varphi}\right)(t) \right| \,f(s,u(s))\, \dif s \\[2pt]
& \le  \int_{\ol A} \, \left| z_p^x(s)- \partial_p\left(\frac{G(\cdot,s)}{\varphi}\right)(t) \right| \,\Phi_{R}(s)\, \dif s.
\end{split}\end{equation*}
Now, from $(C_2)$, we know that for every $\varepsilon\in \bR^+$ there exists $\delta\in\bR^+$ such that $d(x,\kappa(t))<\delta$ implies that
\begin{equation*}
	\left| \Gamma_p\left(G(\cdot,s)\right)(x) - \Gamma_p\left(G(\cdot,s)\right)(\kappa(t)) \right| <	\varepsilon\, w_p(s), \quad \forall\,s\in \ol A,
\end{equation*}
or, which is the same,
\[\left| z_p^x(s)- \partial_p\left(\frac{G(\cdot,s)}{\varphi}\right)(t) \right| <\varepsilon\, w_p(s), \quad \forall\,s\in \ol A. \]
Consequently, since $w_p\,\Phi_R\in \Lsp{1}(\ol A)$, there exists some positive constant $c$ such that
\[\left|{\lim_{\tau\to x}}^\kappa\partial_p\left(\frac{Tu}{\varphi}\right)(\tau)-\partial_p\left(\frac{Tu}{\varphi}\right)(t)\right|< \varepsilon\, c \quad \text{for all } u\in B. \]

Thus, using Theorem~\ref{taae}, we conclude that $T(B)$ is relatively compact in $\cC_{\k,\varphi}^m(\ol A)$ and, therefore,  $T$ is a compact operator.
\end{proof}

\begin{rem}
We must note that in the proof of the previous theorem it is necessary to show that the operator $T$ that we are considering has enough regularity. In particular, we have proved that under hypotheses $(C_1)$--$(C_4)$, operator $T$ maps the space $\cC_{\k,\varphi}^m(\ol A)$ to itself. This means that, given $u\in \cC_{\k,\varphi}^m(\ol A)$, we have proved two different things: first, the regularity of $Tu$, and, second, the asymptotic properties of $Tu$.

We observe that the general hypotheses that we have asked the kernel to satisfy in order to prove the regularity of $Tu$ might be too restrictive or too difficult to check in practice. In this sense, we must take into account the fact that in many examples it will be possible to prove directly the regularity of $Tu$, even if hypotheses $(C_1)$--$(C_4)$ are not satisfied. This will be the case, for example, of integral equations whose origin is a differential equation, as in this case it is clear that the inverse operator of a differential one will always have enough regularity. In fact, this will be the case that we will consider in our example in the last section of this paper.
\end{rem}

Now, following the line of \cite{FigToj}, we will consider an abstract cone in the space $\cC_{\k,\varphi}^m(\ol A)$ defined by
\[K_\alpha=\left\{u\in \cC_{\k,\varphi}^m(\ol A)\ : \ \alpha(u)\ge 0 \right\}, \]
where $\alpha$ is a continuous functional $\alpha\colon \cC_{\k,\varphi}^m(\ol A) \rightarrow \bR$ satisfying the three following properties:
\begin{itemize}
\item[$(P_1)$] $\alpha(u+v)\ge \alpha(u)+\alpha(v)$, for all $u,\, v \in \cC_{\k,\varphi}^m(\ol A)$;
\item[$(P_2)$] $\alpha(\lambda\,u)\ge \lambda \, \alpha(u)$ for all $u\in \cC_{\k,\varphi}^m(\ol A)$, $\lambda\ge 0$;
\item[$(P_3)$] $\left[\alpha(u)\ge 0, \ \alpha(u)\le 0 \right] \Rightarrow u\equiv 0$.
\end{itemize}

In order to choose a cone $K_\alpha$ such that $T$ maps the cone into itself, we will require functional $\alpha$ to satisfy the following condition:
\begin{itemize}
\item[$(C_5)$] For all $u\in K_\alpha$, $Tu\in K_\alpha$.
\end{itemize}


Now we will use the well-know fixed point index theory to prove the existence of a fixed point of operator $T$. In order to do so, we will define some suitable subsets of the cone $K_\alpha$ and give some conditions to ensure that the index of these subsets is either $1$ or $0$.

Let us consider the following subsets
\[K_\alpha^{\beta,\rho}=\{u\in K_\alpha\ : \ \beta(u)<\rho \} \]
and
\[K_\alpha^{\gamma,\rho}=\{u\in K_\alpha\ : \ \gamma(u)<\rho \}, \]
where $\beta$ and $\gamma$ are two continuous functionals $\beta,\, \gamma\colon B \rightarrow \bR$, with $\cC_{\k,\varphi}^m(\ol A) \subset B$ and $\int_{\ol A} |G(t,s)| \dif s \in B$, satisfying the following conditions:
\begin{itemize}
\item[$(C_6)$] For every $u\in K_\alpha$ and $\lambda\in \bR^+$ it holds that $\beta(\lambda\,u)=\lambda\, \beta(u)$ and if $u,\,v\in K_\alpha$ are such that $u\le v$, then $\beta(u)\le \beta(v)$.
\item[$(C_7)$] For every $u,\, v\in K_\alpha$ and $\lambda\in \bR^+$ it holds that $\gamma(u+v)\ge \gamma(u)+\gamma(v)$, $\gamma(\lambda\,u)= \lambda\, \gamma(u)$ and \[\gamma(Tu)\ge \int_{\ol A} \gamma(G(\cdot,s))\, f(s,u(s))\, \dif s.\]
\item[$(C_8)$] $\beta(G(\cdot,s)),\, \gamma(G(\cdot,s)) \in \Lsp{1}(\ol A)$ are positive for every $s\in \ol A$.
\item[$(C_9)$] There exists $e\in K_\alpha \setminus\{0\}$ such that $\gamma(e)\ge 0$.
\item[$(C_{10})$] At least one of the following functions
\begin{equation*}\begin{split}
b\colon \bR^+ &\longrightarrow \bR \\
\rho & \longmapsto 	b(\rho):=\sup\{\beta(u)\ : \ u\in K_\alpha, \ \gamma(u)<\rho  \}
\end{split}\end{equation*}
or
\begin{equation*}\begin{split}
		c\colon \bR^+ &\longrightarrow \bR \\
		\rho & \longmapsto 	c(\rho):=\sup\{\gamma(u)\ : \ u\in K_\alpha, \ \beta(u)<\rho \}
\end{split}\end{equation*}
is well defined, that is, the set on which the supremum is taken is nonempty for every $\rho$ and the supremum is finite.
\end{itemize}

The next lemma compiles some classical results regarding the fixed point index formulated in \cite[Theorems 6.2, 7.3 and 7.11]{GranasDug} in a more general framework.

In particular, given  $X$ a Banach space, $K\subset X$ a cone and $\Omega \subset K$ an arbitrary open subset, $\partial\,\Omega$ will denote the boundary of $\Omega$ in the relative topology in $K$, induced by the topology of $X$. Moreover, let us denote by $\operatorname{Fix}(\cT)$ the set of fixed points of $\cT$.
\begin{lem} \label{l-Amann-fixed-point}
	Let $X$ be a Banach space, $K\subset X$ a cone and $\Omega \subset K$ an arbitrary open subset with $0\in\Omega$. Assume that $\cT\colon \overline \Omega \rightarrow K$ is a compact and compactly fixed operator such that $x\neq \cT x$ for all $x\in\partial\,\Omega$.

	Then the fixed point index $i_K(\cT,\Omega)$ has the following properties:
	\begin{enumerate}
		\item If $x\neq \mu\, \cT x$ for all $x\in\partial\, \Omega$ and for every  $\mu\le 1$, then $i_K(\cT,\Omega)=1.$
		\item If $\Omega$ is bounded and there exists $e\in K\setminus\{0\}$ such that $x\neq \cT x+\lambda\, e$ for all $x\in\partial\, \Omega$ and all $\lambda>0$, then $i_K(\cT,\Omega)=0$.
		\item If $i_K(\cT,\Omega)\neq 0$, then $\cT$ has a fixed point in $\Omega$.
		\item If $\Omega_1$ and $\Omega_2$ are two open and disjoint sets such that $\operatorname{Fix}(\cT) \subset \Omega_1\cup \Omega_2 \subset \Omega$, then
		\[i_K(\cT,\Omega)=i_K(\cT,\Omega_1) + i_K(\cT,\Omega_2). \]
	\end{enumerate}
\end{lem}

\begin{lem}
Assume $T\colon \cC_{\k,\varphi}^m (\ol A) \to \cC_{\k,\varphi}^m (\ol A)$ to be well defined, continuous and compact, that hypotheses $(C_5)$--$(C_6)$ and $(C_8)$ hold and let there exist some $\rho>0$ such that
\begin{equation}\tag{$I_\rho^1$}
	0<f^\rho \beta\left(\int_{\ol A} |G(\cdot,s)|\, \dif s\right) <1,
\end{equation}
where
\[f^\rho=\sup\left\{ \frac{f(t,u(t))}{\rho} \ : \ t\in\ol A, \ u\in K_\alpha, \ \beta(u)=\rho  \right\}. \]
Then $i_{K_\alpha}(T,K_\alpha^{\beta,\rho})=1$.
\end{lem}
\begin{proof}
Let us prove that $Tu\neq \mu\, u $ for all $u\in \partial K_\alpha^{\beta,\rho}$ and every $\mu \ge 1$. Suppose, on the contrary, that there exist some $u\in \partial K_\alpha^{\beta,\rho}$ and some $\mu \ge 1$	such that $\mu\, u(t)=Tu(t)$ for all $t\in \ol A$. In such a case, taking $\beta$ on both sides of the equality and using the fact that
\[Tu(t) =\int_{\ol A}G(t,s)\,f(s,u(s))\, \dif s \le \rho\,f^\rho\, \int_{\ol A} |G(\cdot,s)|\, \dif s, \]
we obtain that
\begin{align*}
\mu \,\rho= & \mu \, \beta(u)= \beta(\mu\,u)= \beta (Tu)=\beta \(\int_{\ol A}G(t,s)\,f(s,u(s))\, \dif s\)\\  = & \rho\, f^\rho\frac{1}{\rho\, f^\rho}\beta \(\int_{\ol A}G(t,s)\,f(s,u(s))\, \dif s\) = \rho\, f^\rho \, \beta \(\int_{\ol A}G(t,s)\frac{f(s,u(s))}{\rho\,f^\rho}\, \dif s\)  \\ \le & \rho \, f^\rho \beta\left(\int_{\ol A} |G(\cdot,s)|\, \dif s\right)<\rho,
\end{align*}
which is a contradiction. Thus, $i_{K_\alpha}(T,K_\alpha^{\beta,\rho})=1$.
\end{proof}

\begin{lem}
Assume $T\colon \cC_{\k,\varphi}^m (\ol A) \to \cC_{\k,\varphi}^m (\ol A)$ to be well defined, continuous and compact, that hypotheses $(C_5)$ and $(C_7)$--$(C_9)$ hold and let there exist some $\rho>0$ such that
\begin{equation}\tag{$I_\rho^0$}
K_\alpha^{\gamma,\rho} \text{ is bounded and } \ 	f_\rho \int_{\ol A} \gamma(G(\cdot,s)) \dif s >1,
\end{equation}
where
\[f_\rho=\inf\left\{ \frac{f(t,u(t))}{\rho} \ : \ t\in \ol A, \ u\in K_\alpha, \ \gamma(u)=\rho  \right\}. \]
Then $i_{K_\alpha}(T,K_\alpha^{\gamma,\rho})=0$.
\end{lem}
\begin{proof}
Let us prove that $u\neq Tu+\lambda\,e $ for all $u\in \partial K_\alpha^{\gamma,\rho}$ and every $\lambda>0$, where $e$ is given im $(C_{10})$. Suppose, on the contrary, that there exist some $u\in \partial K_\alpha^{\gamma,\rho}$ and some $\lambda>0$ such that $u(t)=Tu(t)+\lambda\,e(t)$ for all $t\in \ol A$. In such a case, taking $\gamma$ on both sides of the equality, we obtain that
\begin{equation*}\begin{aligned}
\rho & =\gamma(u)=\gamma(Tu+\lambda\,e) \ge \gamma(Tu)+\lambda\, \gamma(e) \\ & \ge \int_{\ol A} \gamma(G(\cdot,s))\, f(s,u(s))\, \dif s \ge \rho \, f_\rho \int_{\ol A} \gamma(G(\cdot,s))\,  \dif s >\rho,
\end{aligned}\end{equation*}
which is a contradiction. Thus, $i_{K_\alpha}(T,K_\alpha^{\gamma,\rho})=0$.
\end{proof}

Finally, we will give an existence result. We note that, although we will formulate sufficient conditions to ensure the existence of one or two fixed points of operator $T$, it is possible to give similar results to show the existence of three or more fixed points.
\begin{thm}
Let $(C_1)$--$(C_{10})$ hold. Then:
\begin{enumerate}
	\item If the function $b$ given in $(C_{10})$ is well defined and two constants $\rho_1,\,\rho_2\in (0,\infty)$ with $\rho_2>b(\rho_1)$ such that $(I_{\rho_1}^0)$ and $(I_{\rho_2}^1)$ hold, then $T$ has at least a fixed point.
	\item If the function $c$ given in $(C_{10})$ is well defined and two constants $\rho_1,\,\rho_2\in (0,\infty)$ with $\rho_2>c(\rho_1)$ such that $(I_{\rho_1}^1)$ and $(I_{\rho_2}^0)$ hold, then $T$ has at least a fixed point.
	\item If the functions $b$ and $c$ given in $(C_{10})$ are well defined and three constants $\rho_1,\,\rho_2,\, \rho_3\in (0,\infty)$ with $\rho_2>b(\rho_1)$ and $\rho_3>c(\rho_2)$ such that $(I_{\rho_1}^0)$, $(I_{\rho_2}^1)$ and $(I_{\rho_3}^0)$ hold, then $T$ has at least two fixed points.
	\item If the function $b$ and $c$ given in $(C_{10})$ are well defined and three constants $\rho_1,\,\rho_2,\, \rho_3\in (0,\infty)$ with $\rho_2>c(\rho_1)$ and $\rho_3>b(\rho_2)$ such that $(I_{\rho_1}^1)$, $(I_{\rho_2}^0)$ and $(I_{\rho_3}^1)$ hold, then $T$ has at least two fixed points.
\end{enumerate}
\end{thm}
\begin{proof}
The proof follows as a direct consequence of previous lemmas, taking into account that $K_\alpha^{\beta,\rho} \subset K_\alpha^{\gamma,c(\rho)}$ and $K_\alpha^{\gamma,\rho}\subset K_\alpha^{\beta,b(\rho)}$, in case functions $b$ and/or $c$ are well defined.
\end{proof}

\section{An example}

Let $\phi(x,y)= e^{-\frac{1}{2}x^2}$ for every $x,y\in[0,\infty)$. We will consider the following hyperbolic equation:
\begin{equation}\label{eqcal}
	\begin{aligned}
	\frac{\pd^2 u}{\pd x\pd y}(x,y)= & \frac{1}{8}e^{-(x^2+y^2)}+u(x,y)^2,  \quad  (x,y)\in \bR^+\times [0,1].
		\end{aligned}
\end{equation}
 We will look for a positive solution such that there exists
\begin{equation}\label{bc}\lim_{(x,y)\to(\infty,y_0)}\frac{u(x,y)}{\phi(x,y)}\in\bR, \quad \forall\,y_0\in[0,1].\end{equation}

In order to do so, we will consider $\ol A=[0,\infty)\times[0,1]$ the space $\cC_{\k,\varphi}^0(\ol A)$ and $\k:\ol A\to X:=[0,\infty]\times[0,1]$ where $[0,\infty]$ is endowed with the one point compactification topology of $[0,\infty)$.
We will also consider the Green's function
\[G(x,y;t,s)=\begin{dcases} e^{-(x-t)^2}, &  t\in[0,x],\ s\in[0,y],\\ 0 , &  \text{ otherwise.}\end{dcases}\]

Let $f(x,y,u(x,y))=\frac{1}{8} e^{-(x^2+y^2)}+u(x,y)^2$. It is clear that, if we define \[T u(x,y)=\int_{\bR^2}G(x,y;t,s)f(t,s,u(t,s))\dif t\dif s;\quad x,y\in[0,\infty)^2,\] we have that
\[\frac{\pd^2 Tu}{\pd x\pd y}(x,y)=f(x,y,u(x,y)).\]

Let us check that the operator $T:\cC_{\k,\varphi}^0(\ol A)\to\cC_{\k,\varphi}^0(\ol A)$ (defined as the continuous extension of $T$) is well defined, continuous and compact.

\emph{Step 1: Well-definedness:}

If $u\in\cC_{\k,\varphi}^0([0,\infty)^2)$, we have that
\begin{align*}Tu(x,y)= & \int_0^y\int_0^xe^{-(x-t)^2}\(\frac{1}{8}e^{-(t^2+s^2)}+u(t,s)^2\)\dif t\dif s\\= & \frac{\pi  e^{-\frac{x^2}{2}} \operatorname{erf}\left(\frac{x}{\sqrt{2}}\right) \operatorname{erf}(y)}{16 \sqrt{2}}+\int_0^y\int_0^xe^{-(x-t)^2}u(t,s)^2\dif t\dif s,\end{align*}
where $\operatorname{erf}$ denotes the\emph{ error function}
\[{\displaystyle \operatorname {erf} (z)={\frac {2}{\sqrt {\pi }}}\int _{0}^{z}e^{-t^{2}}\,\dif t.}\]
Now, for $y_0\in [0,1]$,
\begin{align*}\lim_{(x,y)\to(\infty,y_0)}\frac{Tu(x,y)}{\phi(x,y)}=& \lim_{(x,y)\to(\infty,y_0)}\left[\frac{\pi  \operatorname{erf}\left(\frac{x}{\sqrt{2}}\right) \operatorname{erf}(y)}{16 \sqrt{2}}+\int_0^y\int_0^xe^{-(x-t)^2}u(t,s)^2\dif t\dif s\right]\\=  & \, \frac{\pi \operatorname{erf}(y_0)}{16 \sqrt{2}}+\lim_{(x,y)\to(\infty,y_0)}\int_0^y\int_0^xe^{-(x-t)^2}u(t,s)^2\dif t\dif s. \end{align*}

Thus, taking into account that $\lim_{(x,y)\to(\infty,y_0)}u(x,y)/\phi(x,y)=L(y_0)\in\bR$, for every $\e\in\bR^+$ there exist $R,\d\in\bR^+$ such that, if $x>R$ and $|y-y_0|<\d$, then $|u(x,y)^2/\phi(x,y)^2-L(y_0)^2|<\e$. Now,
\begin{align*}
&	\int_0^y\int_0^xe^{-(x-t)^2}u(t,s)^2\dif t\dif s \le  	\int_0^{\max\{y_0+\d,1\}}\int_0^\infty e^{-(x-t)^2}u(t,s)^2\dif t\dif s.
\end{align*}
On the other hand,
\begin{align*}
	&	\int_0^y\int_0^xe^{-(x-t)^2}u(t,s)^2\dif t\dif s \ge  	\int_0^{y_0-\d}\int_0^x e^{-(x-t)^2}u(t,s)^2\dif t\dif s\\ = & \int_0^{y_0-\d}\int_0^\infty e^{-(x-t)^2}u(t,s)^2\dif t\dif s-\int_0^{y_0-\d}\int_x^\infty e^{-(x-t)^2}u(t,s)^2\dif t\dif s.
\end{align*}
Since the integrand is bounded, the last integral converges to zero when $R$ tends to infinity (as $x>R$), so we can assume that, if $x>R$ and $|y-y_0|<\d$,
\begin{align*}
	&	\int_0^y\int_0^xe^{-(x-t)^2}u(t,s)^2\dif t\dif s \ge  \int_0^{y_0-\d}\int_0^\infty e^{-(x-t)^2}u(t,s)^2\dif t\dif s-\e.
\end{align*}
Thus, we conclude that
\[\lim_{(x,y)\to(\infty,y_0)}\int_0^y\int_0^xe^{-(x-t)^2}u(t,s)^2\dif t\dif s=\lim_{(x,y)\to(\infty,y_0)}\int_0^{y_0}\int_0^\infty e^{-(x-t)^2}u(t,s)^2\dif t\dif s=0,\]
by the Dominated Convergence Theorem.
Thus,
\begin{align*}\lim_{(x,y)\to(\infty,y_0)}\frac{Tu(x,y)}{\phi(x,y)} = & \frac{\pi \operatorname{erf}(y_0)}{16 \sqrt{2}}.\end{align*}
Hence, $Tu\in\cC_{\k,\varphi}^0(\ol A)$.

\emph{Step 2: Continuity:}

Observe that $(C_3)$ holds, as $f$ is uniformly continuous and, for all $r\in\bR^+$,
for all $y\in\bR$ with $|y|<r$ and a.\,e. $t\in \bR^n$ we have that
\[f(x,y,z\phi(x,y))=\frac{1}{8}e^{-(x^2+y^2)}+z^2e^{-x^2}\le \frac{1}{8}e^{-(x^2+y^2)}+r^2e^{-x^2}:=\Phi_r(x,y),\]
where $\Phi_r\in \Lsp{1}(\ol A)$.

Let $(u_n)_{n\in\bN} \subset \cC_{\k,\varphi}^0 (\ol A)$ be a sequence which converges to $u$ in $\cC_{\k,\varphi}^0 (\ol A)$ and let us show that $(Tu_n)_{n\in\bN}$ converges to $Tu$ in $\cC_{\k,\varphi}^0 (\ol A)$.

The convergence of $(u_n)_{n\in\bN}$ to $u$ in $\cC_{\k,\varphi}^0 (\ol A)$ implies uniform convergence and, from $(C_3)$, $f((x,y),u_n((x,y)))\to f((x,y),u((x,y)))$ uniformly on $(x,y)\in\ol A$.
Take $N\in\bN$ such that for $n\ge N$, $|f((x,y),u_n((x,y)))- f((x,y),u((x,y)))|<\e$. Then, for $n\ge N$,

\begin{align*}
		\left|\frac{Tu_n}{\varphi}(x,y)-\frac{Tu}{\varphi}(x,y) \right|
		& \le  2\e\int_{\ol A}\frac{G(x,y;s,t)}{\varphi(x,y)}\, \dif s=2\e\int_{0}^x\int_{0}^ye^{\frac{1}{2}x^2-(x-t)^2}\, \dif s\dif t\\ & =\e\sqrt{\pi } ye^{\frac{x^2}{2}} \operatorname{erf}(x).
\end{align*}

This way we have proved that $(Tu_n)_{n\in\bN}$ converges to $Tu$ in $\cC_{\k,\varphi}^0 (\ol A)$. Hence, $T$ is a continuous operator.

\emph{Step 3: Compactness:}

Let us consider a bounded set $B\subset \cC_{\k,\varphi}^0 (\ol A)$, that is, such that there exists a positive constant $R$ for which $\|u\|_{\kappa,\varphi}\le R$ for every $u\in B$. We shall prove that $T(B)$ is relatively compact in $\cC_{\k,\varphi}^0 (\ol A)$.

We have that
\begin{equation*}\begin{split}
		\left|\frac{Tu}{\varphi}(x,y) \right| & = \int_{0}^x\int_{0}^y \, \frac{G(x,y;s,t)}{\varphi(x,y)}\, |f(t,s,u(t,s))|\, \dif s\dif t  \le \frac{\sqrt{\pi }}{2}M ye^{\frac{x^2}{2}} \operatorname{erf}(x),
\end{split}\end{equation*}
for every $(x,y)\in \ol A$ and $u\in B$. Therefore, we deduce that $T(B)$ is uniformly bounded.

Now we want to show that for every $(x,y)\in \ol A$ and $\varepsilon\in \bR^+$ there exists some $\delta\in\bR^+$ such that
\[\left|\frac{Tu}{\varphi}(x,y)-\frac{Tu}{\varphi}(t,s)\right|<\varepsilon,\]
if $u\in B$ and $(t,s)\in\ol A$ is such that $\|(x,y)-(t,s)\|<\d$. Observe that
\begin{align*}& \left|\frac{Tu}{\varphi}(x,y)-\frac{Tu}{\varphi}(t,s)\right|\\\le & \left|\frac{\pi \operatorname{erf}\left(\frac{x}{\sqrt{2}}\right) \operatorname{erf}(y)}{16 \sqrt{2}}-\frac{\pi  \operatorname{erf}\left(\frac{t}{\sqrt{2}}\right) \operatorname{erf}(s)}{16 \sqrt{2}}\right|\\ & +\left|\int_0^y\int_0^xe^{\frac{1}{2}x^2-(x-\tau)^2}u(\tau,\sigma)^2\dif \tau \dif \sigma-\int_0^s\int_0^te^{\frac{1}{2}t^2-(t-\tau)^2}u(\tau,\sigma)^2\dif \tau \dif \sigma\right|.\end{align*}
By the linearity of the integral, the continuity of $e^{\frac{1}{2}x^2-(x-\tau)^2}$ and the boundedness of $B$, it is clear that we can bound the previous expression by $\e$ for $(x,y)$ and $(t,s)$ sufficiently close.

On the other hand, we also want to show that, for every $y_0\in [0,1]$ and $\varepsilon\in \bR^+$, there exist some $\delta,M\in\bR^+$ such that
 \[\left|{\lim_{(x,y)\to (\infty,y_0)}}\frac{Tu}{\varphi}(x,y)-\frac{Tu}{\varphi}(t,s)\right|<\varepsilon,\]
 	if $u\in B$ and $(t,s)\in \ol A$ is such that $|s-y_0|<\d$, $t>M$.

In this case, using the same arguments as in Step 1,
\begin{align*}
&\left|\lim_{(x,y)\to  (\infty,y_0)}\frac{Tu}{\varphi}(x,y)-\frac{Tu}{\varphi}(t,s)\right|
\le  \left|\lim_{(x,y)\to (\infty,y_0)}\frac{\pi \operatorname{erf}\left(\frac{x}{\sqrt{2}}\right) \operatorname{erf}(y)}{16 \sqrt{2}}-\frac{\pi  \operatorname{erf}\left(\frac{t}{\sqrt{2}}\right) \operatorname{erf}(s)}{16 \sqrt{2}}\right|
	\\ & +\left|\lim_{(x,y)\to (\infty,y_0)}\int_0^y\int_0^xe^{\frac{1}{2}x^2-(x-\tau)^2}u(\tau,\sigma)^2\dif \tau \dif \sigma-\int_0^s\int_0^te^{\frac{1}{2}t^2-(t-\tau)^2}u(\tau,\sigma)^2\dif \tau \dif \sigma\right|\\ =
	&\left|\frac{\pi \operatorname{erf}(y_0)}{16 \sqrt{2}}-\frac{\pi  \operatorname{erf}\left(\frac{t}{\sqrt{2}}\right) \operatorname{erf}(s)}{16 \sqrt{2}}\right|
 +\left|\int_0^s\int_0^te^{\frac{1}{2}t^2-(t-\tau)^2}u(\tau,\sigma)^2\dif \tau \dif \sigma\right|.
	\end{align*}
For $t$ sufficiently big and $s$ sufficiently close to $y_0$, we have that
\[\left|\int_0^s\int_0^te^{\frac{1}{2}t^2-(t-\tau)^2}u(\tau,\sigma)^2\dif \tau \dif \sigma\right|<\frac{\e}{2},\]
and
\[\left|\frac{\pi \operatorname{erf}(y_0)}{16 \sqrt{2}}-\frac{\pi  \operatorname{erf}\left(\frac{t}{\sqrt{2}}\right) \operatorname{erf}(s)}{16 \sqrt{2}}\right|<\frac{\e}{2},\]
so we get the desired bound.

Thus, using Theorem~\ref{taae}, we conclude that $T(B)$ is relatively compact on $\cC_{\k,\varphi}^0(\ol A)$ and, consequently,  $T$ is a compact operator.

In order to look for a positive solution we consider the cone \[P=\{u\in\cC_{\k,\varphi}^0 (\ol A)\ :\ \a(u)\ge 0\}\] where $\a(u):=\inf u$. $\a$ clearly satisfies the properties $(P_1)-(P_3)$ and it is obvious from the definition of $T$ and the fact that $f,G\ge 0$ that $T$ maps $P$ to $P$.

Now let $\b(u)=\|u\|_\infty$, $\c(u)=0$. Observe that $(C_6)-(C_{10})$ hold. For a given $\rho\in\bR^+$, $\rho\le f^\rho\le (1/8+\rho^2)/\rho$. It holds that
\[\int_{\ol A} |G(x,y;t,s)| \dif t\, \dif s= \frac{1}{2}\, \sqrt{\pi}\,y\,\operatorname{erf}(x) \le \frac{\sqrt{\pi}}{2} <1  \]
and therefore
\begin{equation}
	0<f^\rho \, \beta\left(\int_{\ol A} |G(\cdot,\cdot;t,s)| \dif t\, \dif s\right) < f^\rho  \le (1/8+\rho^2)/\rho .
\end{equation}
Moreover, if $\rho\in \left(\frac{1}{4}(2-\sqrt{2}), \, \frac{1}{4} (2+\sqrt 2)\right)$, then
\[f^\rho \, \beta\left(\int_{\ol A} |G(\cdot,\cdot;t,s)| \dif t\, \dif s\right) <(1/8+\rho^2)/\rho<1,\] so
$i_{K_\alpha}(T,K_\alpha^{\beta,\rho})=1$ and so, by Lemma~\ref{l-Amann-fixed-point}, point 3, there is a solution of problem \eqref{eqcal}--\eqref{bc} with $\|u\|<\rho$.

\section*{Acknowledgements}
The authors would like to acknowledge their gratitude towards Professors Jose Antonio Oubiña and Miguel Domínguez for suggesting some of the bibliography concerning compactifications.

The authors were partially supported by Xunta de Galicia, project ED431C 2019/02, and by the Agencia Estatal de Investigaci\'on (AEI) of Spain under grant MTM2016-75140-P, co-financed by the European Community fund FEDER.
\bibliography{arxiv}
\bibliographystyle{spmpsciper}
\end{document}